\documentclass[a4paper,11pt]{article}
\usepackage{amssymb}
\usepackage{amsmath}
\usepackage{amsfonts}
\usepackage{amsmath}
\usepackage{graphicx}
\numberwithin{equation}{section}
\newtheorem{theorem2}{Theorem}[section]
\newtheorem{lemma2}{Lemma}[section]
\newtheorem{proposition2}{Proposition}[section]

\newtheorem {remark2}{Remark}[section]
\newtheorem {definition2}{Definition}[section]
\newtheorem {theorem3}{Theorem}[section]

\newtheorem {lemma3}{Lemma}[section]

\begin{document}

\title{\textbf{Global Existence and Nonlinear Stability for the Coupled CGL--Burgers Equations for Sequential flames in $\mathbb{R}^N$}}
\author{ Boling Guo\\
Institute of Applied Physics and Computational Mathematics\\
 Beijing, 100088, P. R. China
 \footnote{E-mail: gbl@mail.iapcm.ac.cn}
 \\
 Xinglong Wu\\
Institute of Applied Physics and Computational Mathematics\\
100088, Beijing, China \footnote{E-mail: wxl8758669@yahoo.com.cn}}
\date{}
\maketitle

\begin{abstract}
The present paper is devoted to the study of the global solution and
nonlinear stability to the coupled complex Ginzburg--Landau and
Burgers (CGL--Burgers) equations for sequential flames which
describe the interaction of the excited oscillatory mode and the
damped monotonic mode and are derived from the nonlinear evolution
of the coupled long-scale oscillatory and monotonic instabilities of
a uniformly propagating combustion wave governed by a sequential
chemical reaction, having two flame fronts corresponding to two
reaction zones with a finite separation distance between them. We
first obtain a priori estimation in homogeneous Besov spaces to a
heat equation, thanks to the lemma of  priori estimation, we show
the global solution in a critical Besov space for the Cauchy problem
of Eq.(2.8) if the initial data is small. Next, we obtain the
nonlinear stability by the linearized method as the coefficient
satisfies certain condition. Finally, thanks to two lemmas, we show
 the nonlinear instability for the Cauchy problem of Eq.(3.4), if
we choose certain constants.

\end{abstract}

\section{Introduction}
In this paper, we study the Cauchy problem of the coupled
CGL--Burgers equations for sequential flames \cite{GMBN}
\begin{equation}
\left\{\begin{array}{ll}\partial_{t}P+\nabla Q\cdot\nabla
P-(1+iu)\triangle P=\xi P-(1+iv)|P|^2P-r_1\triangle QP ,
\\\partial_{t}Q+\frac{1}{2}|\nabla Q|^2-m\triangle Q+\kappa|P|^2=0,\; (t,x)\in\mathbb{R}^{+}\times\mathbb{R}^N,
\\P(0,x)=P_0(x), Q(0,x) =Q_{0}(x), \qquad x\in \mathbb{R}^N,\end{array}\right.
\end{equation}
 which describe the interaction of the excited oscillatory mode and the damped monotonic
mode, where $P$ denotes the rescaled complex amplitude of the flame
oscillations, and $Q$ is the deformation of the first front, the
constants
$u=\frac{\delta_i}{\delta_r},v=\frac{\lambda_i}{\lambda_r},
\xi\in\mathbb{R}, m=-\frac{\mu}{\delta_r}$, and
$\kappa=\frac{\eta}{\lambda_r\delta_r}$ are real, otherwise
$r_1=v_1$ is complex. These coefficients as well as the other
parameters in Eq.(1.1) were derived from the original model of
flames governed by a sequential reaction in \cite{GMBN}. Applying
$\nabla$  to Eq.$(1.1)_2$, letting $\Omega=\nabla Q$. Substitute
$\Omega$ into Eq.$(1.1)_1$, we obtain
\begin{equation}
\left\{\begin{array}{ll}\partial_{t}P+\Omega\cdot\nabla
P-(1+iu)\triangle P=\xi P-(1+iv)|P|^2P-r_1Pdiv\Omega  ,
\\\partial_{t}\Omega+\Omega\cdot\nabla\Omega-m\triangle\Omega+\kappa\nabla(|P|^2)=0,\; (t,x)\in\mathbb{R}^{+}\times\mathbb{R}^N,
\\P(0,x)=P_0(x), \Omega(0,x) =\Omega_{0}(x), \qquad x\in \mathbb{R}^N,\end{array}\right.
\end{equation}
where $\Omega\cdot\nabla
\Omega=\sum_{i=1}^N\Omega_i\partial_{x_i}\Omega$,
$div\Omega=\sum_{i=1}^N\partial_{x_i}\Omega_i$.
\par
As well known, uniformly propagating planar premixed flame fronts
can become unstable leading to a lot of spatial or spatio-temporal
flame structures, such as cellular flames, pulsating flames,
spinning and spiral flames, flames with traveling or standing waves
on the flame fronts and many others can be found in \cite{Wi,ZBLM}
as well as the references cited therein. Premixed gaseous flames
generally involve many reactions and reactants. In premixed gaseous
flames with a one-stage chemical reaction (the case of a single
limiting reactant which  reacts to form products), there exist two
kinds of diffusional thermal instability \cite{Bu}, which occur
depending on the value of Lewis number $L=\kappa_T/D$, where
$\kappa_T$ denotes thermal diffusivity of the reactive mixture, and
$D$ is the diffusivity of the deficient reacting component. If
$L<L_{c1}<1$, the flame front is unstable to monotonic long-scale
perturbations that develop into cellular structures whose
spatio-temporal evolution is governed by the Kuramoto--Sivashunsky
(KS) equation \cite{Si}. KS equation describes stationary spatially
periodic patterns as well as  further transitions to increasingly
complex spatio-temporal evolution, eventually leading to chaotic or
weakly turbulent behavior. If $1<L_{c2}<L$, an oscillatory
instability occurs with a nonzero wavenumber and frequency at the
instability threshold, and its nonlinear evolution is described by
the complex Ginzburg-Landau (CGL) equations \cite{OM}.
\par
However, the chemical kinetics governing the structure of flames can
be quite complex and  can lead to qualitative flames behavior, which
can not be accounted for by a one-stage kinetic model. The effect of
complex chemical kinetics on the dynamics of propagating flame
fronts is governed by the two-stage sequential reaction, which leads
to the occurrence of two flame fronts with various regimes of
propagation depending on the heat releases, the reaction rate
constants, and the activation energies of the chemical reactions
\cite{BR,KL,MM,PL}. One of the possible regimes of propagation is
that of two planar uniformly propagating fronts with a constant
separation distance between them. A linear stability analysis of
this regime for the sequential reaction problem proved that, as in
the case of a one-stage reaction, the flame fronts can become
unstable with respect to either monotonic or oscillatory
instabilities \cite{Pe}. The nonlinear interaction between the
monotonic and oscillatory modes of instability was investigated in
reaction diffusion systems \cite{Ki,PWD}, in Rayleigh--Benard and
Benard--Marangoni convection \cite{CGLL,FR}. If the monotonic mode
was excited  with a nonzero wavenumber, the resulting nonlinear
evolution was described by coupled  real and CGL equations.
\par
The nonlinear interaction between the monotonic and oscillatory
modes of instability of the two uniformly propagating flame fronts
can be governed by the following coupled KS-CGL equations
\cite{GMBN}
\begin{equation}
\left\{\begin{array}{ll}\partial_{t}P+\nabla Q\cdot\nabla
P-(1+iu)\triangle P=\xi P-(1+iv)|P|^2P-(r_1\triangle
Q+r_2g\triangle^2Q)P ,
\\\partial_{t}Q+\frac{1}{2}|\nabla Q|^2+\kappa|P|^2=m\triangle Q-g\triangle^2Q,\; (t,x)\in\mathbb{R}^{+}\times\mathbb{R}^N,
\\P(0,x)=P_0(x), Q(0,x) =Q_{0}(x), \qquad x\in \mathbb{R}^N,\end{array}\right.
\end{equation}
which describes the interaction of the evolving monotonic and
oscillatory modes of instability of the uniform flame fronts. where
$r_2=v_2\delta_r,\;g=\frac{|\gamma_r(l_A-l_A^0)|\chi}{\delta_r^2}$,
and other coefficients present as Eq.(1.1). The oscillatory and
monotonic modes are excited as $m<0$ and $\xi=1$. In the case $m<0$
and $\xi=-1$, the oscillatory mode is damped and the monotonic mode
is excited. If $m>0$ and $\xi=1$, the oscillatory mode is excited
and the monotonic mode is damped, in this case of letting $g=0$,
then Eq.(1.3) becomes the coupled CGL--Burgers equations. In
contrast to the coupled KS-CGL equations, the coupled CGL--Burgers
equations which describes the interaction between the excited
oscillatory mode and the damped monotonic mode is an asymptotically
valid description of the system behavior.  Eq.(1.1) is  CGL equation
that describes the weakly nonlinear evolution of a long-scale
instability if no coupled monotonic mode $Q$ \cite{AK,AKW,CM}, while
Eq.$(1.2)_2$ becomes the Burgers equation if we let the coefficient
$\kappa=0$. There are some works concerning the coupled CGL--Burgers
equations. In \cite{GMBN}, for the 1D version of Eq.(1.1), they
obtain that the coupling with the Goldstone mode leads to new types
of dynamics of a planar pulsating flame resulting from the Hopf
bifurcation at zero wavenumber by the method of numerical
computation. In \cite{GNM}, for 2D of Eq.(1.1), the authors focus on
the effect of the coupling between the two modes, $P$ and $Q$, on
the evolution of traveling waves, they show that the coupled system
exhibits new types of instabilities as well as new dynamical
behavior, including bound states of two or four spirals, ¡°liquid
spiral¡± states, superspiral structures, oscillating cellular
structures separated by chaotically merging and splitting domain
walls.
\par
 The remainder of this paper is organized as follows. In Section 2, we first recall
 the Littlewood--Paley decomposition, the Bony decomposition,
the definition and properties of Besov space. Next, we obtain a
priori estimation in homogeneous Besov spaces to Eq.(2.1).  Thanks
to these lemmas, we establish global solution for the Cauchy problem
of Eq.(2.8) in  Besov spaces by Theorem 2.1. By virtue of some
remarks, we get the global solution to Eq.(1.1). In Section 3, we
first recall the definition of spectral stability. Next, we obtain
the nonlinear stability by linearization to Eq.(3.1). In Section 4,
by virtue of two lemmas, we get the the nonlinear instability to the
Cauchy problem of Eq.(3.4) if we choose certain constants. \\

Notation: Let $\widetilde{L}^\rho_T(\dot{B}_{p,r}^\sigma)$ denote
the set of functions $u$ such that
$$\|u\|_{\widetilde{L}^\rho_T(\dot{B}_{p,r}^\sigma)}:=\left\|\left(2^{k\sigma}\|\dot{\triangle}_ku\|_{L_T^\rho(L^p)}\right)_{k\in\mathbb{Z}}\right\|_{l^r},$$
 for $\sigma\in\mathbb{R},T>0$, and $(p,r,\rho)\in[1,\infty]^3$. Then
 via the Minkowski inequality, we have
 $$\|u\|_{\widetilde{L}^\rho_T(\dot{B}_{p,r}^\sigma)}\leq C\|u\|_{L^\rho_T(\dot{B}_{p,r}^\sigma)}\qquad \text{if}\quad 1\leq\rho\leq r,$$
otherwise,
$$\|u\|_{L^\rho_T(\dot{B}_{p,r}^\sigma)}\leq C\|u\|_{\widetilde{L}^\rho_T(\dot{B}_{p,r}^\sigma)}\qquad \text{if}\quad 1\leq r\leq\rho $$
If for any $T>0$, $u\in\widetilde{L}^\rho_T(\dot{B}_{p,r}^\sigma)$,
then we have $u\in\widetilde{L}^\rho(\dot{B}_{p,r}^\sigma)$, where
$$\|u\|_{\widetilde{L}^\rho(\dot{B}_{p,r}^\sigma)}:=\left(\sum_{k\in\mathbb{Z}}(2^{k\sigma}\|\dot{\triangle}_ku
\|_{L^\rho(\mathbb{R}^+;L^p)})^r\right)^{1/r}.$$

\section{Global existence of the solution}
In this section, in order to establish the global solution of the
Cauchy problem for  Eq.(1.1) in Besov spaces. First, for the
convenience of the readers, we recall some facts on
the Littlewood--Paley decomposition, the Bony decomposition and some useful lemmas.\\
\begin{proposition2} \cite{B-C-D} There exists a couple
of $C^{\infty}$ functions $(\chi, \varphi)$ valued in [0, 1], such
that $\chi$ is supported in the ball $\mathcal {B}=\{\xi\in
\mathbb{R}^N;|\xi|\leq\frac{4}{3}\}$, and $\varphi$ is supported in
the annulus $\mathcal {C}=\{\xi\in
\mathbb{R}^N;\frac{3}{4}\leq|\xi|\leq\frac{8}{3}\}$. Moreover,
$$\chi(\xi)+\sum_{q\in\mathbb{N}}\varphi(2^{-q}\xi)=1,\qquad \forall\xi\in\mathbb{R}^N,$$
$$q\geq1\Longrightarrow\;\text{supp}\;\chi(\cdot)\cap\text{supp}\;\varphi(2^{-q}\cdot)=\emptyset,$$
$$\text{supp}\;\varphi(2^{-q}\cdot)\cap\text{supp}\;\varphi(2^{-p}\cdot)=\emptyset,\;\text{if}
\;|p-q|\geq2$$
and$$\frac{1}{3}\leq\chi^2(\xi)+\sum_{q\in\mathbb{N}}\varphi^2(2^{-q}\xi)\leq1,\;\forall\xi\in\mathbb{R}^N.$$
\end{proposition2}
Let $\tilde{h}=\mathcal {F}^{-1}\chi$ and
$h=\mathcal{F}^{-1}\varphi$. Then the nonhomogeneous dyadic blocks
$\Delta_q$ and the nonhomogeneous low-off operators $S_q$ can be
defined as follows:
$$\Delta_{-1}u=S_0u\quad\text{and}\quad \Delta_qu=0,\quad \text{if}\quad q\leq-2,$$
$$\Delta_qu=\varphi(2^{-q}D)u=2^{qN}\int_{\mathbb{R}^N}h(2^qy)u(x-y)dy,\;\text{if}\;q\geq0,$$
$$S_qu=\sum_{p\geq-1}^{q-1}\Delta_pu=\chi(2^{-q}D)u=2^{qN}\int_{\mathbb{R}^N}\tilde{h}(2^qy)u(x-y)dy.$$
Moreover, if $u,v\in\mathcal{S}'(\mathbb{R}^N)$, then we have
$$\Delta_p\Delta_q=0\quad \text{if}\quad |p-q|\geq2,$$
$$\Delta_q(S_{p-1}u\Delta_pv)=0\quad \text{if}\quad |p-q|\geq5.$$
Furthermore, for all $u\in \mathcal{S}'(\mathbb{R}^N)$, one can
easily check that
$$u=\sum_{q\in\mathbb{Z}}\Delta_qu \quad\text{in}\quad \mathcal{S}'(\mathbb{R}^N).$$
The homogeneous dyadic blocks $\dot{\Delta}_q$ and the homogeneous
low-off operators $\dot{S}_q$ are defined for all $q\in\mathbb{Z}$
by
$$\dot{\Delta}_qu=\varphi(2^{-q}D)u=2^{qN}\int_{\mathbb{R}^N}h(2^qy)u(x-y)dy,$$
$$\dot{S}_qu=\sum_{p\leq q-1}\dot{\Delta}_pu=\chi(2^{-q}D)u=2^{qN}\int_{\mathbb{R}^N}\tilde{h}(2^qy)u(x-y)dy.$$
Then the nonhomogeneous Besov space $B_{p,r}^s(\mathbb{R}^N)$ and
the homogeneous Besov space $\dot{B}_{p,r}^s(\mathbb{R}^N)$ is
defined as follows:
\begin{definition2} Let $s\in\mathbb{R},\;p,r\in[1,\infty]$, we set
$$B_{p,r}^s=\left\{u\in\mathcal{S}^{'}(\mathbb{R}^N);\|u\|_{B_{p,r}^s}=
\left(\sum_{k=-1}^{\infty}2^{ksr}\|\Delta_{k}u\|_{L^p}^r\right)^{1/r}<\infty\right\},$$
where $\Delta_k$ are the nonhomogeneous dyadic blocks. If
$s=\infty,\;B_{p,r}^\infty=\cap_{\sigma\in\mathbb{R}}B_{p,r}^\sigma$.
$$\dot{B}_{p,r}^s=\left\{u\in\mathcal{S}^{'}(\mathbb{R}^N);\|u\|_{\dot{B}_{p,r}^s}=
\left\|(\|2^{ks}\Delta_{k}u\|_{L^p})_{k\in\mathbb{Z}}\right\|_{l^r(\mathbb{Z})}<\infty\right\},$$
where $\dot{\Delta}_k$ is the homogeneous dyadic blocks.
 If
$s=\infty,\;\dot{B}_{p,r}^\infty=\cap_{\sigma\in\mathbb{R}}\dot{B}_{p,r}^\sigma$.

\end{definition2}
For $u,v\in\mathcal{S}'(\mathbb{R}^N)$, we have the Bony
decomposition as follows:
\begin{definition2} Let $u,v\in\mathcal{S}'(\mathbb{R}^N)$. Denote
$$T_{u}v=\sum_{q\geq-1}\sum_{p\geq-1}^{q-2}\Delta_pu\Delta_qv=\sum_{q\geq-1}S_{q-1}u\Delta_qv$$
and$$R(u,v)=\sum_{q\geq-1}\Delta_qu\widetilde{\Delta_q}v\;\text{with}\;\widetilde{\Delta_q}=
\Delta_{q-1}+\Delta_q+\Delta_{q+1}.$$Then formally, we have the
nonhomogeneous Bony decomposition
$$uv=T_{u}v+T_{v}u+R(u,v).$$
Similarly,  the homogeneous Bony decomposition gives rise by
$$uv=\dot{T}_{u}v+\dot{T}_{v}u+\dot{R}(u,v),$$ where $\dot{T}_{u}v=\sum_{q\in\mathbb{Z}}\sum_{p\leq
q-2}\dot{\Delta}_pu\dot{\Delta}_qv=\sum_{q\in\mathbb{Z}}\dot{S}_{q-1}u\dot{\Delta}_qv$
and
$$\dot{R}(u,v)=\sum_{q\in\mathbb{Z}}\dot{\Delta}_qu\widetilde{\dot{\Delta}_q}v\;\text{with}\;\widetilde{\dot{\Delta}_q}=
\dot{\Delta}_{q-1}+\dot{\Delta}_q+\dot{\Delta}_{q+1}.$$
\end{definition2}
 We now state the result concerning continuity of the inhomogeneous paraproduct operator $T$ and
 the remainder operator $R$.
 \begin{lemma2} \cite{B-C-D} There exists a constant $C$ such that
 for any couple of  real numbers $(s,t)$ with $t$ negative  and any $(p,r,r_1,r_2)\in[1,\infty]^4$, we
 have, for any $(u,v)\in L^\infty\times B_{p,r}^s$,
 $$\|T_uv\|_{B_{p,r}^s}\leq C^{1+s}\|u\|_{L^\infty}\|v\|_{B_{p,r}^s},$$
 for any $(u,v)\in B_{\infty,r_1}^{t}\times B_{p,r_2}^s$ and
 $\frac{1}{r}:=\min\{1,\frac{1}{r_1}+\frac{1}{r_2}\}$
 $$\|T_uv\|_{B_{p,r}^{s+t}}\leq \frac{C^{1+|s+t|}}{-t}\|u\|_{B_{\infty,r_1}^t}\|v\|_{B_{p,r_2}^s}.$$
 Moreover, assume $(s_1,s_2)$ be real, and
 $(p_1,p_2,r_1,r_2)\in[1,\infty]^4$ such that$$\frac{1}{p}:=+\frac{1}{p_1}+\frac{1}{p_2}\leq1\quad\text{and}\;
 \frac{1}{r}:=\frac{1}{r_1}+\frac{1}{r_2}\leq1.$$ If $s_1+s_2>0$,
 then we imply, for any $(u,v)\in B_{p_1,r_1}^{s_1}\times
 B_{p_2,r_2}^{s_2}$, $$\|R(u,v)\|_{B_{p,r}^{s_1+s_2}}\leq \frac{C^{1+|s_1+s_2|}}{s_1+s_2}\|u\|_{B_{p_1,r_1}^{s_1}}\|v\|_{B_{p_2,r_2}^{s_2}}.$$
If $r=1$ and $s_1+s_2\geq0$, we have, for any $(u,v)\in
B_{p_1,r_1}^{s_1}\times B_{p_2,r_2}^{s_2}$,
$$\|R(u,v)\|_{B_{p,\infty}^{s_1+s_2}}\leq C^{1+|s_1+s_2|}\|u\|_{B_{p_1,r_1}^{s_1}}\|v\|_{B_{p_2,r_2}^{s_2}}.$$
\end{lemma2}
\begin{remark2} Similar to the case of inhomogeneous, for the
homogeneous paraproduct operator $\dot{T}$ and  the remainder
operator $\dot{R}$, we can get the same result.
\end{remark2}

\begin{lemma2} \cite{B-C-D} The following properties hold:\\
i) Density: if $p,r<\infty$, then $\mathcal{S}(\mathbb{R}^N)$ is
dense in $B_{p,r}^s$. Let the space $\mathcal {S}_0(\mathbb{R}^N)$
denotes the function in $\mathcal{S}(\mathbb{R}^N)$ whose Fourier
transforms are supported away from $0$. Then the space $\mathcal {S}_0(\mathbb{R}^N)$ is dense in $\dot{B}_{p,r}^s$.\\
ii) Generalized derivatives: let $f\in C^\infty(\mathbb{R}^N)$ be a
homogeneous function of degree $m\in\mathbb{R}$ away from a
neighborhood of the origin. There exists a constant $C$ depending
only on $f$ and such that $\|f(D)u\|_{B_{p,r}^s}\leq
C\|u\|_{B_{p,r}^{s+m}}.$\\ iii) Sobolev embeddings: if $p_1\leq p_2$
and $r_1\leq r_2$, then $B_{p_1,r_1}^s\hookrightarrow
B_{p_2,r_2}^{s-N(\frac{1}{p_1}-\frac{1}{p_2})}$. If $s_1<s_2,1\leq
p\leq\infty$ and $1\leq r_1,r_2\leq\infty$, then the embedding
$B_{p,r_2}^{s_2}\hookrightarrow B_{p,r_1}^{s_1}$ is locally
compact.\\
iv) Algebraic properties: for $s>0,\; B_{p,r}^s\cap L^\infty$ is an
algebra. Moreover $(B_{p,r}^s$ is an algebra$)\Leftrightarrow$
$(B_{p,r}^s\hookrightarrow L^\infty)$
$\Leftrightarrow (s>N/p$ or $(s\geq N/p$ and $r=1))$.\\
 v) Fatou property: if $(u_n)_{n\in\mathbb{N}}$ is a bounded sequence of
 $B_{p,r}^s$ which tends to $u$ in $\mathcal{S}^{'}$, then $u\in B_{p,r}^s$
 and $\|u\|_{B_{p,r}^s}\leq\lim_{n\rightarrow\infty}\inf\|u_n\|_{B_{p,r}^s}$.\\
vi) Complex interpolation: if $u\in B_{p,r}^s\cap B_{p,r}^{s_1}$ and
$\theta\in[0,1],\;p,r\in[1,\infty]$, then $u\in B_{p,r}^{\theta
s+(1-\theta)s_1}$ and
$$\|u\|_{B_{p,r}^{\theta s+(1-\theta)s_1}}\leq\|u\|_{B_{p,r}^s}^\theta\|u\|_{B_{p,r}^{s_1}}^{1-\theta}.$$
vii) Let $m\in\mathbb{R}$ and $f$ be a $S^m$-multiplier, i.e. $f :
\mathbb{R}^n \mapsto \mathbb{R}$ is smooth and satisfies that for
all multi-index $\alpha$, there exists a constant $C_{\alpha}$ such
that $\forall \xi\in\mathbb{R}^n\;, |\partial^{\alpha}f(\xi)|\leq
C_\alpha(1+|\xi|)^{m-|\alpha|}.$ Then for all $s \in\mathbb{R}$ and
 $p,r\in[1,\infty],$ the operator $f(D)$ is continuous from $B_{p,r}^s$ to
 $B_{p,r}^{s-m}$.
\end{lemma2}
\begin{remark2} Properties ii), v), vi), vii) hold for the
homogeneous spaces $\dot{B}_{p,r}^s$, and the following properties
hold for the homogeneous spaces; Sobolev embeddings: if $p_1\leq
p_2$ and $r_1\leq r_2$, then $\dot{B}_{p_1,r_1}^s\hookrightarrow
\dot{B}_{p_2,r_2}^{s-N(\frac{1}{p_1}-\frac{1}{p_2})}$.  Moreover,
for $2\leq p<\infty$, $\dot{B}_{p,2}^0\hookrightarrow L^p$.
Algebraic properties: for $s>0,\;\dot{ B}_{p,r}^s\cap L^\infty$ is
an algebra. Moreover $(\dot{B}_{p,r}^s$ is an
algebra$)\Leftrightarrow$ $(\dot{B}_{p,r}^s\hookrightarrow
L^\infty)$ $\Leftrightarrow (s=N/p$ and $r=1))$.
\end{remark2}

By the Bony decomposition, we can infer the following estimate.
\begin{lemma2}\cite{B-C-D} For any positive real number $s$ and
any $p,r\in[1,\infty]$, the space $L^\infty\cap B_{p,r}^{s}$ is an
algebra, and there exists a constant $C$ such that
$$\|uv\|_{B_{p,r}^s}\leq \frac{C^{s+1}}{s}(\|u\|_{L^\infty}\|v\|_{B_{p,r}^s}+\|v\|_{L^\infty}\|u\|_{B_{p,r}^s}).$$
\end{lemma2}
\begin{lemma2}\cite{B-C-D} Assume $f$ is a  smooth function such that $f(0)=0$, $s>0$ and
 $p,r\in[1,\infty]$. If $u$ belongs to the space $L^\infty\cap
B_{p,r}^{s}$, then we have
$$\|f(u)\|_{B_{p,r}^s}\leq
C(s,f',\|u\|_{L^\infty})\|u\|_{B_{p,r}^s}.$$Moreover, if $f$ belongs
to $C_b^\infty(\mathbb{R})$ and $u$ belongs to
$B_{\infty,\infty}^{-1}$, then we also obtain
$$\|f(u)\|_{B_{p,r}^s}\leq
C(s,f,\|\nabla u\|_{B_{\infty,\infty}^{-1}})\|u\|_{B_{p,r}^s}.$$
\end{lemma2}
 Finally, taking
advantage of $$f(u)-f(v)=(u-v)\int_0^1f'(v+\theta(u-v))d\theta,$$ we
can get the following lemma.
\begin{lemma2}\cite{B-C-D} Assume $f$ is a  smooth function such that $f(0)=0$, $s>0$ and
 $p,r\in[1,\infty]$. For any couple $(u,v)$ belongs to the space
$L^\infty\cap B_{p,r}^{s}$,  we have the function $f(u)-f(v)$
belongs to $L^\infty\cap B_{p,r}^{s}$ and
\begin{equation*}\begin{split}\|f(u)-f(v)\|_{B_{p,r}^s}&\leq
C(s,f'',\|u\|_{L^\infty},\|v\|_{L^\infty})(\sup_{\theta\in[0,1]}
\|u+\theta(v-u)\|_{L^\infty}\\&\|u-v\|_{B_{p,r}^s}+\|u-v\|_{L^\infty}\sup_{\theta\in[0,1]}
\|u+\theta(v-u)\|_{B_{p,r}^s}).
\end{split}\end{equation*}
\end{lemma2}
\par
Consider the following $n$-dimension linear equation
\begin{equation}
\left\{\begin{array}{ll}\partial_tf-\mu(1+iu)\triangle f=g
,&(t,x)\in\mathbb{R}^{+}\times\mathbb{R}^N,
\\ f(0,x) =f_{0}(x),&x\in \mathbb{R}^N,\end{array}\right.
\end{equation}
where $f(t,x)=(f_1,\cdots,f_n)$ is the vector field, the external
source term $g=g(t,x)$ and the initial $f_0$ are known data. The
diffusion constant $\mu$ is  positive. Applying the partial Fourier
transformation with respect  to the space  variable,  we have
\begin{equation}
\hat{f}(t,\xi)=e^{-\mu(1+iu)|\xi|^2t}\hat{f}_0(\xi)+\int_0^te^{-\mu(1+iu)|\xi|^2(t-\tau)}\hat{g}(\tau,\xi)d\tau.
 \end{equation}
If we introduce the semi-group $\{e^{t\Delta}\}_{t\geq0}$, then
(2.2) is equivalent to
\begin{equation}
f(t,x)=e^{\mu(1+iu)t\triangle}f_0(x)+\int_0^te^{\mu(1+iu)(t-\tau)\triangle}g(\tau,x)d\tau.
 \end{equation}
Analogy to the proof of Lemma 2.4 on page 54 in \cite{B-C-D}, we can
obtain the following result.
\begin{lemma2} Assume $p\in[1,\infty]$, $\mathcal {C}$ be an annulus, and supp$\hat{f}$ be the support
set of $f$. If supp$\hat{f}\subset\lambda\mathcal {C}$, for some
$\lambda>0$. Then there exists constants $c$ and $C$, for any
positive number $t,\mu$, we obtain
\begin{equation*}\|e^{\mu(1+iu)t\triangle}f\|_{L^p}\leq Ce^{-c\mu
t\lambda^2}\|f\|_{L^p}.
\end{equation*}
\end{lemma2}
By Lemma 2.6, we can get an important lemma to Eq.(2.1).

\begin{lemma2} Given $(\rho,p,r)\in[1,\infty]^3$ and
$\sigma\in\mathbb{R}$, Assume that the initial data $f_0\in
\dot{B}_{p,r}^\sigma$ and the external source term $g\in
\widetilde{L}^\rho_T(\dot{B}_{p,r}^{\sigma-2+\frac{2}{\rho}})$, for
a fixed $T>0$. Then for all $\rho_1\in [\rho,\infty]$, Eq.(2.1) has
a unique solution
$f\in\widetilde{L}^{\rho_1}_T(\dot{B}_{p,r}^{\sigma+\frac{2}{\rho_1}})\cap
\widetilde{L}^\infty_T(\dot{B}_{p,r}^{\sigma})$ with the initial
data $f_0$. Moreover, there exists a constant $C$ such that
\begin{equation}\mu^{\frac{1}{\rho}}\|f\|_{\widetilde{L}^{\rho_1}_T(\dot{B}_{p,r}^{\sigma+\frac{2}{\rho_1}})}\leq
C\|f_0\|_{\dot{B}_{p,r}^\sigma}
+C\mu^{\frac{1}{\rho}-1}\|g\|_{\widetilde{L}^\rho_T(\dot{B}_{p,r}^{\sigma-2+\frac{2}{\rho}})}.
\end{equation}
In addition, if $r$ is finite, then the solution $f\in \mathcal
{C}([0,T];\dot{B}_{p,r}^{\sigma})$.
\end{lemma2}
\textit{Proof.} Applying the dyadic blocks $\dot{\triangle}_k$ to
Eq.(2.1) to yield
\begin{equation}
\left\{\begin{array}{ll}\partial_t\dot{\triangle}_kf-\mu(1+iu)\triangle
\dot{\triangle}_kf=\dot{\triangle}_kg
,&(t,x)\in\mathbb{R}^{+}\times\mathbb{R}^N,
\\ \dot{\triangle}_kf(0,x) =\dot{\triangle}_kf_{0}(x),&x\in \mathbb{R}^N,
\end{array}\right.
\end{equation}
Using formula (2.3) to Eq.(2.5) to obtain
\begin{equation}
\dot{\triangle}_kf(t,x)=e^{\mu(1+iu)t\triangle}\dot{\triangle}_kf_0(x)+\int_0^te^{\mu(1+iu)(t-\tau)\triangle}\dot{\triangle}_kg(\tau,x)d\tau.
\end{equation}
Taking advantage of Lemma 2.6 to (2.6), for some $c>0$, we have
\begin{equation}\begin{split}
\|\dot{\triangle}_kf&\|_{L^{\rho_1}_T(L^p)}\leq
C\|e^{-c\mu2^{2k}t}\|\dot{\triangle}_kf_0\|_{L^p}\|_{L_T^{\rho_1}}+\|e^{-c\mu2^{2k}t}*\|\dot{\triangle}_kg(t)\|_{L^p}\|_{L_T^{\rho_1}}\\&
\leq
C\left(\frac{1-e^{-c\mu\rho_12^{2k}T}}{c\mu\rho_12^{2k}}\right)^{\frac{1}{\rho_1}}\|\dot{\triangle}_kf_0\|_{L^p}+
C\left(\frac{1-e^{-c\mu\rho_22^{2k}T}}{c\mu\rho_22^{2k}}\right)^{\frac{1}{\rho_2}}\|\dot{\triangle}_kg\|_{L_T^\rho(L^p)},
\end{split}\end{equation}
where we used Young's inequality and
$1+\frac{1}{\rho_1}=\frac{1}{\rho}+\frac{1}{\rho_2}$.\\
Therefore, we get (2.4) by taking the $l^r({\mathbb{Z}})$ norm to
(2.7). Moreover, one can easily get the
 solution $f\in \mathcal {C}([0,T];\dot{B}_{p,r}^{\sigma})$, if
$r$ is finite.$\qquad \Box$\\
\par
 Next, we consider the following modified CGL--Burgers equations
\begin{equation}
\left\{\begin{array}{ll}(\partial_{t}-(1+iu)\triangle)
P=-\Omega\cdot\nabla P-(1+iv)|P|^2P-r_1Pdiv\Omega +f_1 ,
\\(\partial_{t}-m\triangle)\Omega=-\Omega\cdot\nabla\Omega-\kappa\nabla(|P|^2)+f_2,\; (t,x)\in\mathbb{R}^{+}\times\mathbb{R}^N,
\\P(0,x)=P_0(x), \Omega(0,x) =\Omega_{0}(x), \qquad x\in \mathbb{R}^N,\end{array}\right.
\end{equation}
where $m>0$ and other coefficients present as Eq.(1.1), $f_1$ and
$f_2$ are the external source terms.
\par
Now, we present the main result of this paper.
\begin{theorem2}
  Assume $N\geq2,p\in[1,2N[$, and $(P_0,\Omega_0)\in (\dot{B}_{p,1}^{N/p-1})^2$. Then there exists a constant $s$ for all external force
$f_i \in\widetilde{L}^{1}(\dot{B}_{p,1}^{N/p-1}),\;i=1,2$ such that
$$\|(P_0,\Omega_0)\|_{\dot{B}_{p,1}^{N/p-1}}+\sum_{i=1}^2\|f_i\|_{\widetilde{L}^{1}(\dot{B}_{p,1}^{N/p-1})}\leq s,$$ system
(2.8) has  a unique global solution $(P,\Omega)$ belongs to
$\widetilde{L}^{\infty}(\dot{B}_{p,1}^{N/p-1})\cap\widetilde{L}^{1}(\dot{B}_{p,1}^{N/p+1})$
with the initial data $(P_0,\Omega_0)$ which satisfies
$$\|(P,\Omega)\|_{\widetilde{L}^{\infty}(\dot{B}_{p,1}^{N/p-1})\cap\widetilde{L}^{1}(\dot{B}_{p,1}^{N/p+1})}\leq 2s.$$
Moreover, the solution $(P,\Omega)$ belongs to $\mathcal
{C}(\mathbb{R}^+;\dot{B}_{p,1}^{N/p-1})$.
\end{theorem2}
\begin{remark2} The spaces $\widetilde{L}^{\infty}(\dot{B}_{p,1}^{N/p-1})\cap\widetilde{L}^{1}(\dot{B}_{p,1}^{N/p+1})$
are the scaling invariance in the following transformation. If
$(P(t,x),\Omega(t,x))$ is the solution to Eq.(2.8) with the initial
$(P_0(x),\Omega_0(x))$ and external force term $f_i$, then for any
$\lambda>0$, $P_\lambda(t,x)=\lambda P(\lambda^2t,\lambda x)$ and
$\Omega_\lambda(t,x)=\lambda \Omega(\lambda^2t,\lambda x)$ solve
Eq.(2.8) with initial data
$(P_{0,\lambda},\Omega_{0,\lambda})=(\lambda P_0(\lambda x),\lambda
\Omega_0(\lambda x))$ and
$f_{i,\lambda}=\lambda^3f_i(t\lambda^2,\lambda x)$. Moreover, we
have
$$\|(P_\lambda,\Omega_\lambda)\|_{\widetilde{L}^{\infty}(\dot{B}_{p,1}^{N/p-1})\cap\widetilde{L}^{1}(\dot{B}_{p,1}^{N/p+1})}=
\|(P,\Omega)\|_{\widetilde{L}^{\infty}(\dot{B}_{p,1}^{N/p-1})\cap\widetilde{L}^{1}(\dot{B}_{p,1}^{N/p+1})}$$
and
$$\|f_{i,\lambda}\|_{\widetilde{L}^{1}(\dot{B}_{p,1}^{N/p-1})}=\|f_i\|_{\widetilde{L}^{1}(\dot{B}_{p,1}^{N/p-1})}.$$
\end{remark2}
\begin{remark2} If we choose a compact domain $K$ of $\mathbb{R}^N$.
Let $\dot{B}_{p,r}^\sigma(K)$ (resp., $B_{p,r}^\sigma(K)$) denote
the set of distributions $f$ in $\dot{B}_{p,r}^\sigma(\mathbb{R}^N)$
(resp., $B_{p,r}^\sigma(K)$), the support of which is included in
$K$. By Proposition 2.93 on page 108 in \cite{B-C-D}, if $p<N$, then
the spaces $\dot{B}_{p,1}^{N/p-1}(K)$ and $B_{p,1}^{N/p-1}(K)$
coincide. Furthermore, we have
$$\|\xi(P_1-P_2)\|_{\widetilde{L}^{1}(\dot{B}_{p,1}^{N/p-1})}\lesssim
\|P_1-P_2\|_{\widetilde{L}^{1}(\dot{B}_{p,1}^{N/p+1})}.$$ Therefore,
if $m>0$, we can get the same result of Theorem 2.1 in the spaces
$$\widetilde{L}^{\infty}(\dot{B}_{p,1}^{N/p-1}(K))\bigcap\widetilde{L}^{1}(\dot{B}_{p,1}^{N/p+1}(K))$$
to Eq.(1.2).
\end{remark2}
\begin{remark2}Due to the estimation (2.4) is true for
the morn $L^{\rho}_{T}(\dot{B}_{p,r}^{s})$. Therefore, if we take
place
$\widetilde{L}^{\infty}(\dot{B}_{p,1}^{N/p-1}),\;\widetilde{L}^{1}(\dot{B}_{p,1}^{N/p+1})$
by $L^{\infty}(\dot{B}_{p,1}^{N/p-1}),L^{1}(\dot{B}_{p,1}^{N/p+1})$
in Theorem 2.1, respectively, the result is also right. Moreover,
$$\widetilde{L}^{\infty}(\dot{B}_{p,1}^{N/p-1})\cap\widetilde{L}^{1}(\dot{B}_{p,1}^{N/p+1})\hookrightarrow
L^{\infty}(\dot{B}_{p,1}^{N/p-1})\cap
L^{1}(\dot{B}_{p,1}^{N/p+1}).$$
\end{remark2}
\par
 We introduce the solution $(P_{F},\Omega_F)$ of the
following linear equation
\begin{equation}\left\{\begin{array}{ll}\partial_tP_F-(1+iu)\triangle
P_F=f_1 ,&(t,x)\in\mathbb{R}^{+}\times\mathbb{R}^N,
\\ \partial_t\Omega_F-m\triangle \Omega_F=f_2,
\\(P_F, \Omega_F)|_{t=0}=(P_0, \Omega_0), &x\in
\mathbb{R}^N.\end{array}\right.
\end{equation}
By the virtue of  Duhamel's formula, $(P,\Omega)$ is a solution of
Eq.(2.8) if and only if $\left(\begin{array}{cc} P  \\\Omega
\end{array}\right)= \left(\begin{array}{cc}
P_{F}+F(P)\\
\Omega_F+F(\Omega)
\end{array}\right)$
with
\begin{equation}\left\{\begin{array}{ll}(\partial_t-(1+iu)\triangle)F(P)=-\Omega\cdot\nabla P
-(1+iv)|P|^2P+r_1P div\Omega ,
\\ (\partial_t-m\triangle)F(\Omega)=-\Omega\cdot\nabla\Omega-\kappa\nabla(|P|^2),\;(t,x)\in\mathbb{R}^{+}\times\mathbb{R}^N,
\\(F(P), F(\Omega))|_{t=0}=(0,0), \;x\in
\mathbb{R}^N.\end{array}\right.
\end{equation}
The existence of solution to Eq.(2.8) in Theorem 2.1 is a
straightforward corollary of the following result.
\begin{lemma2}Given $\beta\geq1,N\geq2$ and $p\in[1,2N[$, define the set $X_{p}^{s}(\beta)$ of the  functions $u$ such
that for any $T>0$,
 $$u\in
 \widetilde{L}^{\infty}_T(\dot{B}_{p,1}^{N/p-1})\bigcap\widetilde{L}^{1}_T(\dot{B}_{p,1}^{N/p+1}),\quad
 \text{and}\quad
 \|u\|_{\widetilde{L}^{\infty}(\dot{B}_{p,1}^{N/p-1})\cap\widetilde{L}^{1}(\dot{B}_{p,1}^{N/p+1})}\leq\beta s$$
 if the initial data $u_0$ satisfies $\|u_0\|_{\dot{B}_{p,1}^{N/p-1}}\leq
 s$, for some small enough $s$. Then the mapping  $\left(\begin{array}{cc} P  \\\Omega
\end{array}\right)\longmapsto \left(\begin{array}{cc}
P_{F}+F(P)\\
\Omega_F+F(\Omega)
\end{array}\right)$ maps $X_{p}^{s}(\beta)$ into
 $X_{p}^{s}(\beta)$, if the initial data $(P_0,\Omega_0)\in (X_{p}^{s}(\beta))^2$. More precisely, for any solutions $(P_1,\Omega_1)$ and $(P_2,\Omega_2)$ in
 $X_{p}^{s}(\beta)$, the following result holds
 \begin{equation}\begin{split}\|(F(P_1)-F(P_2),&F(\Omega_1)-F(\Omega_2))\|_{\widetilde{L}^{\infty}(\dot{B}_{p,1}^{N/p-1})\cap\widetilde{L}^{1}(\dot{B}_{p,1}^{N/p+1})}
 \\&\leq\frac{1}{2}\|(P_1-P_2,\Omega_1-\Omega_2)\|_{\widetilde{L}^{\infty}(\dot{B}_{p,1}^{N/p-1})\cap\widetilde{L}^{1}(\dot{B}_{p,1}^{N/p+1})}.
 \end{split}\end{equation}
\end{lemma2}
\textit{Proof.} Since $(P_1,\Omega_1)$ and $(P_2,\Omega_2)$ are two
solutions to Eq.(2.8), it is obviously  that
$(F(P_1)-F(P_2),F(\Omega_1)-F(\Omega_2))$ solves the following
equation
\begin{equation}\left\{\begin{array}{ll}(\partial_t-(1+iu)\triangle)(F(P_1)-F(P_2))=Q(P_1,P_2,\Omega_1,\Omega_2),
\\ (\partial_t-m\triangle)(F(\Omega_1)-F(\Omega_2))=W(P_1,P_2,\Omega_1,\Omega_2),
\\(F(P_1)-F(P_2),F(\Omega_1)-F(\Omega_2))|_{t=0}=(0,0),\end{array}\right.
\end{equation}
where
$Q(P_1,P_2,\Omega_1,\Omega_2)=-(1+iv)(|P_1|^2P_1-|P_2|^2P_2)-(\Omega_1\cdot
\nabla P_1-\Omega_2\cdot\nabla P_2)+r_1(div\Omega_1 P_1-div
\Omega_2P_2)$ and
$W(P_1,P_2,\Omega_1,\Omega_2)=-(\Omega_1\cdot\nabla\Omega_1-\Omega_2\cdot\nabla\Omega_2)-\kappa\nabla(|P_1|^2-|P_2|^2)$.
\par
Taking advantage of Lemma 2.7 to Eq.$(2.12)_1$ with
$\sigma=\frac{N}{p}-1,\;\rho=1$ and $r=1$, we obtain
\begin{equation}\begin{split}\|(F(P_1)-F(P_2))\|_{\widetilde{L}^{\rho_1}(\dot{B}_{p,1}^{N/p-1+2/\rho_1})}&\leq
c\|Q\|_{\widetilde{L}^{1}(\dot{B}_{p,1}^{N/p-1})}\\&\lesssim(I+II+III)
\end{split}\end{equation} with
$I=\||P_1|^2P_1-|P_2|^2P_2\|_{\widetilde{L}^{1}(\dot{B}_{p,1}^{N/p-1})}$,
$II=\|\Omega_1\cdot \nabla P_1-\Omega_2\cdot\nabla
P_2\|_{\widetilde{L}^{1}(\dot{B}_{p,1}^{N/p-1})}$ and
$III=\|div\Omega_1 P_1-div
\Omega_2P_2\|_{\widetilde{L}^{1}(\dot{B}_{p,1}^{N/p-1})}$. \par We
now estimate $I,II,III$ respectively. Due to
$$|P_1|^2P_1-|P_2|^2P_2=|P_1|^2(P_1-P_2)+P_1P_2(\bar{P_1}-\bar{P_2})+|P_2|^2(P_1-P_2)$$
and the Bony decomposition
$$|P_1|^2(P_1-P_2)=\dot{T}_{|P_1|^2}(P_1-P_2)+\dot{T}_{(P_1-P_2)}|P_1|^2+\dot{R}(|P_1|^2,(P_1-P_2)).$$
Thanks to Lemma 2.1, we have
\begin{equation*}\begin{split}\|\dot{T}_{|P_1|^2}(P_1-P_2)\|_{\dot{B}_{p,1}^{N/p-1}}&\leq\|P_1\|_{L^\infty}^2\|P_1-P_2\|_{\dot{B}_{p,1}^{N/p-1}},
\\&\leq\|P_1\|_{\dot{B}_{p,1}^{N/p}}^2\|P_1-P_2\|_{\dot{B}_{p,1}^{N/p-1}},\end{split}\end{equation*}

\begin{equation*}\begin{split}\|\dot{T}_{(P_1-P_2)}|P_1|^2\|_{\dot{B}_{p,1}^{N/p-1}}&\leq\|P_1-P_2\|_{\dot{B}_{\infty,\infty}^{-1}}
\||P_1|^2\|_{\dot{B}_{p,1}^{N/p}}\\&\leq\|P_1\|_{\dot{B}_{p,1}^{N/p}}^2\|P_1-P_2\|_{\dot{B}_{p,1}^{N/p-1}},
\end{split}\end{equation*}
$$\|\dot{R}(|P_1|^2,(P_1-P_2))\|_{\dot{B}_{p,1}^{N/p-1}}\leq
\|P_1\|_{\dot{B}_{p,1}^{N/p}}^2\|P_1-P_2\|_{\dot{B}_{p,1}^{N/p-1}},\;
\text{if}\;1\leq p<N,$$
\begin{equation*}\begin{split}
\|\dot{R}(|P_1|^2,(P_1-P_2))\|_{\dot{B}_{p,1}^{N/p-1}}&\leq\|\dot{R}(|P_1|^2,(P_1-P_2))\|_{\dot{B}_{p/2,1}^{2N/p-1}}\\&\leq\|P_1\|_{\dot{B}_{p,1}^{N/p}}^2\|P_1-P_2\|_{\dot{B}_{p,1}^{N/p-1}},\;
\text{if}\;2\leq p<2N.
 \end{split}\end{equation*}
 Therefore, by virtue of the above inequalities, applying interpolation theorem, it follows that
 \begin{equation}
 \||P_1|^2(P_1-P_2)\|_{\dot{B}_{p,1}^{N/p-1}}\leq\|P_1\|_{\dot{B}_{p,1}^{N/p-1}}\|P_1\|_{\dot{B}_{p,1}^{N/p+1}}\|P_1-P_2\|_{\dot{B}_{p,1}^{N/p-1}}.
\end{equation}
Similarly, along the same lines to estimate above inequality, we
immediately imply that
\begin{equation}\left\{\begin{array}{ll}\|P_1P_2(\bar{P_1}-\bar{P_2})\|_{\dot{B}_{p,1}^{N/p-1}}\leq\|P_1\|_{\dot{B}_{p,1}^{N/p}}
\|P_2\|_{\dot{B}_{p,1}^{N/p}} \|P_1-P_2\|_{\dot{B}_{p,1}^{N/p-1}},
\\ \||P_2|^2(P_1-P_2)\|_{\dot{B}_{p,1}^{N/p-1}}\leq
\|P_2\|_{\dot{B}_{p,1}^{N/p}}^2 \|P_1-P_2\|_{\dot{B}_{p,1}^{N/p-1}}.
\end{array}\right.
\end{equation}
Therefore, combining (2.14) with (2.15), as $N\geq2$, $p\in[1,2N[$,
we deduce
 \begin{equation}I\leq\left(\sum_{j=1}^2\|P_j\|_{\widetilde{L}^{\infty}(\dot{B}_{p,1}^{N/p-1})}\|P_j\|_{\widetilde{L}^{1}(\dot{B}_{p,1}^{N/p+1})}\right)
 \|P_1-P_2\|_{\widetilde{L}^{\infty}(\dot{B}_{p,1}^{N/p-1})}.
\end{equation}
Note that $\Omega_1\cdot \nabla P_1-\Omega_2\cdot\nabla
P_2=\Omega_1\cdot \nabla( P_1-P_2)+(\Omega_1-\Omega_2)\cdot \nabla
P_2 $ and $div\Omega_1 P_1-div\Omega_2
P_2=P_1div(\Omega_1-\Omega_2)+div\Omega_2(P_1-P_2).$
\par
 Now, argument as we did in proving (2.14), we end up with
\begin{equation}II\leq\|\Omega_1-\Omega_2\|_{\widetilde{L}^{\infty}(\dot{B}_{p,1}^{N/p-1})}\|P_2\|_{\widetilde{L}^{1}(\dot{B}_{p,1}^{N/p+1})}+
\|\Omega_1\|_{\widetilde{L}^{\infty}(\dot{B}_{p,1}^{N/p-1})}\|P_1-P_2\|_{\widetilde{L}^{1}(\dot{B}_{p,1}^{N/p+1})},
\end{equation}
\begin{equation}III\leq\|P_1\|_{\widetilde{L}^{\infty}(\dot{B}_{p,1}^{N/p-1})}\|\Omega_1-\Omega_2\|_{\widetilde{L}^{1}(\dot{B}_{p,1}^{N/p+1})}+
\|\Omega_2\|_{\widetilde{L}^{1}(\dot{B}_{p,1}^{N/p+1})}\|P_1-P_2\|_{\widetilde{L}^{\infty}(\dot{B}_{p,1}^{N/p-1})}.
\end{equation}
Combining (2.13), (2.16), (2.17) with (2.18), for $\rho_1=1$ and
$\rho_1=\infty$, one can easily check that
\begin{equation}\begin{split}\|F(P_1)-F(P_2)&\|_{X_1\cap X_2}\leq
\|\Omega_1\|_{X_1}\|P_1-P_2\|_{X_2}\\&+\|\Omega_1-\Omega_2\|_{X_1}\|P_2\|_{X_2}+\|P_1\|_{X_1}\|\Omega_1-\Omega_2\|_{X_2}\\
&+\left(\sum_{j=1}^2\|P_j\|_{X_1}\|P_j\|_{X_2}+\|\Omega_2\|_{X_2}\right)
\|P_1-P_2\|_{X_1},
\end{split}\end{equation}
where $X_1=\widetilde{L}^{\infty}(\dot{B}_{p,1}^{N/p-1})$ and
$X_2=\widetilde{L}^{1}(\dot{B}_{p,1}^{N/p+1})$.
\par
Taking advantage of Lemma 2.7 to Eq.$(2.12)_2$ with
$\sigma=\frac{N}{p}-1,\;\rho=1$ and $r=1$, we obtain
\begin{equation}\|(F(P_1)-F(P_2))\|_{\widetilde{L}^{\rho_1}(\dot{B}_{p,1}^{N/p-1+2/\rho_1})}\leq
c\|W\|_{\widetilde{L}^{1}(\dot{B}_{p,1}^{N/p-1})}\leq c(IV+V),
\end{equation}
where
$IV=\|\Omega_1\cdot\nabla\Omega_1-\Omega_2\cdot\nabla\Omega_2\|_{\widetilde{L}^{1}(\dot{B}_{p,1}^{N/p-1})}$,
and $V=\|\nabla(
|P_1|^2-|P_2|^2)\|_{\widetilde{L}^{1}(\dot{B}_{p,1}^{N/p-1})}$.
Along the same lines to estimate (2.14), we immediately imply that
\begin{equation}IV\leq
\|\Omega_1\|_{X_1}\|\Omega_1-\Omega_2\|_{X_2}+\|\Omega_1-\Omega_2\|_{X_1}\|\Omega_2\|_{X_2}.
\end{equation}
Since $|P_1|^2-|P_2|^2=P_1(\bar{P_1}-\bar{P_2})+(P_1-P_2)\bar{P_2}$,
and $\dot{B}_{p,1}^{N/p}\hookrightarrow L^\infty$ is an algebra.
Applying interpolation theorem and Cauchy--Schwarz inequality, we
have
\begin{equation}\begin{split}\|\nabla(
|P_1|^2-|P_2|^2)&\|_{\dot{B}_{p,1}^{N/p-1}}\leq\||P_1|^2-|P_2|^2\|_{\dot{B}_{p,1}^{N/p}}\\&\leq
(\|P_1\|_{\dot{B}_{p,1}^{N/p}}+\|P_2\|_{\dot{B}_{p,1}^{N/p}})\|P_1-P_2\|_{\dot{B}_{p,1}^{N/p}}\\&\leq
\left(\sum_{j=1}^2\|P_i\|_{\dot{B}_{p,1}^{N/p-1}}\right)\|P_1-P_2\|_{\dot{B}_{p,1}^{N/p+1}}\\&+
\left(\sum_{j=1}^2\|P_i\|_{\dot{B}_{p,1}^{N/p+1}}\right)\|P_1-P_2\|_{\dot{B}_{p,1}^{N/p-1}}.
\end{split}\end{equation}
Thus, we can get the estimation of $V$ as follows
\begin{equation}V\leq
(\|P_1\|_{X_1}+\|P_2\|_{X_1})\|P_1-P_2\|_{X_2}+(\|P_1\|_{X_2}+\|P_2\|_{X_2})\|P_1-P_2\|_{X_1},
\end{equation}
where $X_1=\widetilde{L}^{\infty}(\dot{B}_{p,1}^{N/p-1})$ and
$X_2=\widetilde{L}^{1}(\dot{B}_{p,1}^{N/p+1})$.
\par
Plugging (2.21), (2.23) into (2.20), taking $\rho_1=1$and
$\rho_1=\infty$, we imply that
\begin{equation}\begin{split}\|F(\Omega_1)-F(\Omega_2)&\|_{X_1\cap X_2}\lesssim
\|\Omega_1\|_{X_1}\|\Omega_1-\Omega_2\|_{X_2}+\left(\sum_{j=1}^2\|P_j\|_{X_1}\right)\|P_1-P_2\|_{X_2}\\&+
\|\Omega_2\|_{X_2}\|\Omega_1-\Omega_2\|_{X_1}+\left(\sum_{j=1}^2\|P_j\|_{X_2}\right)\|P_1-P_2\|_{X_1}.
\end{split}\end{equation}
Adding (2.19) with (2.24) to obtain
\begin{equation}\begin{split}\|(F(P_1)-&F(P_2),F(\Omega_1)-F(\Omega_2))\|_{X_1\cap
X_2}\\&\leq c\left(\sum_{j=1}^2(\|P_j\|_{X_1\cap
X_2}+\|\Omega_j\|_{X_1\cap
X_2}+\|P_j\|_{X_1}\|P_j\|_{X_2})\right)\\&\times
\|(P_1-P_2,\Omega_1-\Omega_2)\|_{X_1\cap X_2},
\end{split}\end{equation}
where $\|(A,B)\|_{X_1\cap
X_2}=(\|A\|_{X_1}+\|A\|_{X_2})+(\|B\|_{X_1}+\|B\|_{X_2})$.
\par
By the assumption of Lemma 2.8, one can easily get
\begin{equation*}\begin{split} c(\sum_{j=1}^2(\|P_j\|_{X_1\cap
X_2}+\|\Omega_j\|_{X_1\cap X_2}+&\|P_j\|_{X_1}\|P_j\|_{X_2}))
\\&\leq2c\beta s(2+\beta s)=\alpha<1,
\end{split}\end{equation*}
if we let $s$  small enough. The above inequality implies that
\begin{equation*}\begin{split}\|(F(P_1)-&F(P_2),F(\Omega_1)-F(\Omega_2))\|_{X_1\cap
X_2}\leq \alpha\|(P_1-P_2,\Omega_1-\Omega_2)\|_{X_1\cap X_2}.
\end{split}\end{equation*}
\par
 Next, we shall prove the existence of solution to Eq.(2.8), i.e. if
 $(P,\Omega)\in (X_p^s(\beta))^2$, so does $(P_F+F(P),\Omega_F+F(\Omega))\in
 (X_p^s(\beta))^2$.
 \par On one hand, by Lemma 2.7 to $(2.9)_1$ and $(2.10)_1$, we have
 \begin{equation*}\left\{\begin{array}{ll}\|P_F\|_{X_1\cap X_2}\leq
 c(\|P_0\|_{X_1}+\|f_1\|_{\widetilde{L}^1(B_{p,1}^{N/p-1})};\\
\|F(P)\|_{X_1\cap X_2}\leq
c(\|\Omega\|_{X_1}\|P\|_{X_2}+\|\|_{X_1}^2\|P\|_{X_2}).
 \end{array}\right.\end{equation*}
 Therefore, without loss generality, we let $f_1=0$. Then
 \begin{equation}\begin{split}\|P_F+F(P)\|_{X_1\cap X_2}&\leq
 c(\|P_0\|_{X_1}+\|\Omega\|_{X_1}\|P\|_{X_2}+\|P\|_{X_1}^2\|P\|_{X_2})\\&\leq
 c\left(s+(\beta s)^2+(\beta s)^3\right)\\&\leq\beta s,
 \end{split}\end{equation}
 where we let $s$ small enough, and $\beta$ is a fix real number.\\
On the other hand, by Lemma 2.7 to $(2.9)_2$ and $(2.10)_2$, we
imply
\begin{equation*}\left\{\begin{array}{ll}\|\Omega_F\|_{X_1\cap X_2}\leq
 c\|\Omega_0\|_{X_1}+\|f_2\|_{\widetilde{L}^1(B_{p,1}^{N/p-1})};\\
\|F(\Omega)\|_{X_1\cap X_2}\leq
c(\|\Omega\|_{X_1}\|\Omega\|_{X_2}+\|P\|_{X_1}\|P\|_{X_2}).
 \end{array}\right.\end{equation*}
 Therefore, with no loss generality, we let $f_2=0$. Then
 \begin{equation}\begin{split}\|\Omega_F+F(\Omega)\|_{X_1\cap X_2}&\leq
 c(\|\Omega_0\|_{X_1}+\|\Omega\|_{X_1}\|\Omega\|_{X_2}+\|P\|_{X_1}\|P\|_{X_2})\\&\leq
 c\left(s+2(\beta s)^2\right)\\&\leq\beta s,
 \end{split}\end{equation}
 where we let $s$ small enough, and $\beta$ is a fix real number.
  This completes the proof of Lemma 2.8. $\qquad\Box$\\
 \textit{Proof of Theorem 2.1.} Thanks to Lemma 2.8, we can get the
 existence of solution to Eq.(2.8) by the contraction mapping argument, we will complete the proof of
 Theorem 2.1, if we prove the uniqueness. Indeed, let
 $(P_1,\Omega_1)$, $(P_2,\Omega_2)$ be two solutions to Eq.(2.8)
 with initial $(P_{1,0},\Omega_{1,0}),\; (P_{2,0},\Omega_{2,0})$, respectively.
 Then, similar the process of (2.26) and (2.27), by virtue of (2.11), we obtain
 $$(P_1-P_2,\Omega_1-\Omega_2)=(P_{1F}-P_{2F},\Omega_{1F}-\Omega_{2F})+(F(P_1)-F(P_2),F(\Omega_1)-F(\Omega_2))$$
 satisfies the following inequality
 \begin{equation}\begin{split}\|(P_1-P_2,\Omega_1-\Omega_2)&\|_{X_1\cap
 X_2}\leq\|(P_{1F}-P_{2F},\Omega_{1F}-\Omega_{2F})\|_{X_1\cap X_2}\\&+\|(F(P_1)-F(P_2),F(\Omega_1)-F(\Omega_2))\|_{X_1\cap
 X_2}\\&\leq C\|(P_{1,0}-P_{2,0},\Omega_{1,0}-\Omega_{2,0})\|_{X_1\cap
 X_2}\\&+\frac{1}{2}\|(P_1-P_2,\Omega_1-\Omega_2)\|_{X_1\cap X_2},
 \end{split}\end{equation}
where $X_1=\widetilde{L}^{\infty}(\dot{B}_{p,1}^{N/p-1})$ and
$X_2=\widetilde{L}^{1}(\dot{B}_{p,1}^{N/p+1})$. One can easily check
the uniqueness from the term (2.28).  This completes the proof of Theorem 2.1.$\qquad\Box$\\

\section{Nonlinear Stability by Linearization}
In this section, we will establish the the nonlinear stability  of
the Cauchy problem for Eq.(1.1) by linearization. First, for the
convenience of the readers, we recall an useful definition.\\
\begin{definition2} The plane wave is said to be spectrally stable
if the spectrum is bounded in the complex plane by the curve
$Re\lambda=-C(Im\lambda)^2$, for some $C>0$.
\end{definition2}
\par
 Next, letting $P=r(t,x)e^{i\vartheta(t,x)}$, we can obtain from (1.2) that the following CGL--Burgers equations
\begin{equation}
\left\{\begin{array}{ll}r_t+\Omega\cdot\nabla r=\triangle
r-u(2\nabla r\cdot\nabla
\vartheta+r\triangle\vartheta)+r(1-r^2-|\nabla\vartheta|^2-s_1div\Omega),
\\r\vartheta_t+r\Omega\cdot\nabla\vartheta=u(\triangle r-r|\nabla\vartheta|^2)+2\nabla r\cdot\nabla\vartheta+r(\triangle\vartheta+vr^2-s_2div\Omega),
\\\Omega_t=m\triangle \Omega-\Omega\cdot\nabla \Omega-\kappa\nabla (r^2),(t,x)\in\mathbb{R}^+\times\mathbb{T}\end{array}\right.
\end{equation}
where the constant $m>0$, and the other coefficients are the
function of $r$.
\par
Consider the plane wave solution of the CGL--Burgers equations of
one-dimension, which is given by
\begin{equation*}
\left\{\begin{array}{ll}P=r_0e^{i\vartheta_0x},
\\\Omega=w_0.\end{array}\right.
\end{equation*}
Substituting it into (3.1), by simply calculation, we imply that
\begin{equation*}
\left\{\begin{array}{ll}r_0^2+\vartheta_0^2=1,
\\u\vartheta_0^2+vr_0^2=0.\end{array}\right.
\end{equation*}
If we consider the following  perturbation of the wave which is
assumed to satisfy
\begin{equation}
\left\{\begin{array}{ll}r=r_0+\rho,
\\\vartheta=\vartheta_0x+\phi,
\\\Omega=w_0+h.\end{array}\right.
\end{equation}
Then, we have the following result.
\begin{theorem3}If the domain $\mathbb{T}\in \mathbb{R}$ is bounded. Assume $m>0$. Then there exists constants $r_0,\vartheta_0$
and functions $u,v$ such that the plane wave is spectrally stable.
Moreover, letting $\pi=(\rho,\phi,h)$, if
$\|\pi_0\|_{L^1}+\|\pi_0\|_{H^{s+1}},\;s>\frac{1}{2}$ is
sufficiently small and $u=0$, then the perturbations satisfy
$$\|\pi\|_{H^{s+1}}\leq C(1+t)^{-\frac{1}{2}(\frac{3}{2}+s)}\|\pi_0\|_{H^{s+1}}.$$
\end{theorem3}
We break the argument into several lemmas.
\begin{lemma2} Suppose $m>0$. Then there exists constants $r_0,\vartheta_0$
such that the plane wave is spectrally stable.
\end{lemma2}

\textit{Proof.} Plugging (3.2) into (3.1), letting $\kappa=0$, it
follows that
\begin{equation}
\left\{\begin{array}{ll}\rho_t=\rho_{xx}-w_0\rho_xx-h\rho_x-u[2\rho_x(\vartheta_0+\phi_x)+
(r_0+\rho)\phi_{xx}]+\\
\qquad (r_0+\rho)[1-(r_0+\rho)^2- (\vartheta_0+\phi_x)^2-s_1h_x],
\\(r_0+\rho)\phi_t+(r_0+\rho)(w_0+h)(\vartheta_0+\phi_x)=2\rho_x(\vartheta_0+\phi_x)+\\
\qquad
u[\rho_{xx}-(r_0+\rho)(\vartheta_0+\phi_x)^2]+(r_0+\rho)[\phi_{xx}-v(r_0+\rho)^2-s_2h_x],
\\h_t=mh_{xx}-(w_0+h)h_x.\end{array}\right.
\end{equation}
Assume $u(r)=c_0+c_1r$. By the formulation of Taylor expansion
\begin{equation*}
\left\{\begin{array}{ll}m(r)=m(r_0)+m'(r_0)\rho+m''(r_0)\rho^2/2+\mathcal
{O}(\rho^3),
\\s_1(r)=s_1(r_0)+s_1'(r_0)\rho+s_1''(r_0)\rho^2/2+\mathcal
{O}(\rho^3),
\\s_2(r)=s_2(r_0)+s_2'(r_0)\rho+s_2''(r_0)\rho^2/2+\mathcal
{O}(\rho^3).\end{array}\right.
\end{equation*}
Therefore, inserting (3.3) into (3.2), by simply computation, it
follows that
\begin{equation}
\left\{\begin{array}{ll}\rho_t=\rho_{xx}-(c_0+c_1r_0)r_0\phi_{xx}-[w_0+2\vartheta_0(c_0+c_1r_0)]\rho_x-2\vartheta_0r_0\phi_x\\
\qquad-r_0s_1(r_0)h_x-2r_0^2\rho-2\vartheta_0c_1\rho\rho_x-2u\phi_x\rho_x-(c_0+c_1r_0)\rho
\phi_{xx}-h\rho_{x}\\
\qquad+r_0[(s_1(r_0)-s_1(r))h_x-\rho^2-\phi_x^2]-\rho[2r_0\rho+\rho^2+2\vartheta_0\phi_x+\phi_x^2+s_1h_x],
\\\phi_t=c_1\rho_{xx}+\phi_{xx}-[2\vartheta_0(c_0+c_1r_0)+w_0]\phi_x-[c_1\vartheta_0^2+r_0^2v'(r_0)+2r_0v(r_0)]\rho\\
\qquad-s_2(r_0)h_x-\vartheta_0h-h\phi_x-w_0\vartheta_0
-(c_0+c_1r_0)(\vartheta_0^2+\phi_x^2)\\\qquad+\frac{c_0}{r_0+\rho}\rho_{xx}-
c_1(2\vartheta_0\phi_x+\phi_x^2)\rho-\frac{2\rho_x}{r_0+\rho}(\vartheta_0+\phi_x)+[s_2(r_0)-s_2(r)]h_x\\
\qquad-v(r)\rho^2-r_0^2[v(r)-v'(r_0)\rho]-2r_0\rho[v(r_0)-v(r)],
\\h_t=mh_{xx}-w_0h_x-hh_x.\end{array}\right.
\end{equation}
Note that (3.4) gives the following linearized equation
\begin{equation}\pi_t=\mathcal{A}\pi_{xx}+\mathcal{B}\pi_x+\mathcal{C}\pi,
\end{equation}
where the vector $\pi$, and the $3\times3$ matrices $\mathcal
{A},\mathcal {B},\mathcal {C}$ are
$$\pi=\left(\begin{array}{ccc}
\rho \\
\phi \\
 h
\end{array}\right),
\mathcal {A}=\left(\begin{array}{ccc}

1&-r_0(c_0+c_1r_0)&0\\
c_1&1&0 \\
0&0&m
\end{array}\right),$$
$$\mathcal {B}=\left(\begin{array}{ccc}
-w_0-2\vartheta_0(c_0+c_1r_0)&-2\vartheta_0r_0&-s_1(r_0)r_0\\
0&-2\vartheta_0(c_0+c_1r_0)-w_0&-s_2(r_0) \\
0&0&-w_0
\end{array}\right),$$
$$\mathcal {C}=\left(\begin{array}{ccc}
-2r_0^2&0&0\\
-[c_1\vartheta_0^2+r_0^2v'(r_0)+2r_0v(r_0)]&0&-\vartheta_0 \\
0&0&0
\end{array}\right).$$
\par

Define the operator $\mathcal {L}$ as follows
$$\mathcal {L}=\mathcal {A}\partial_x^2+\mathcal {B}\partial_x+\mathcal {C}.$$
 As is well known \cite{He}, the spectrum of $\mathcal {L}$ in any
 $L^p$ space, $p\in[1,\infty]$, is bounded by the curve
 \begin{equation}
C_\lambda=\{\lambda;|-k^2\mathcal {A}+ik\mathcal {B}+(\mathcal
{C}-\lambda I_3)|=0,k\in\mathbb{T}\}.
 \end{equation}
A lengthy calculation shows that for $\lambda\in C_\lambda$, such
that
$$\lambda_1(k)=-k^2m-w_0ki,$$
\begin{equation}\begin{split}&\lambda_{2,3}(k)=-[r_0^2+k^2+(w_0+2\vartheta_0(c_0+c_1r_0))ki]\pm\\&\sqrt{r_0^4-c_1r_0(c_0+c_1r_0)k^4+r_0[c_0\vartheta_0^2+r_0^2v'(r_0)+2r_0v(r_0)]k+
2r_0\vartheta_0c_1k^3i},\end{split}\end{equation}

from which one can easily check that
\begin{equation}\begin{split}
Re(\lambda_{2,3}(k))=-(r_0^2+k^2)\pm\frac{[\sqrt{a^2+b^2}+a]^{1/2}}{\sqrt{2}},
\end{split}\end{equation}
where
$a=r_0^4-c_1r_0(c_0+c_1r_0)k^4+r_0[c_0\vartheta_0^2+r_0^2v'(r_0)+2r_0v(r_0)]k$,
and $b=2r_0\vartheta_0c_1k^3.$
\par
If $k^2m>0$, $2(r_0^2+k^2)^2\geq a$, and
$4(r_0^2+k^2)^4-4a(r_0^2+k^2)^2>b^2$, then we deduce that
$$Re\lambda_1(k)<0,\quad\text{and}\quad Re\lambda_{2,3}(k)<0.$$
Therefore, the plane wave of the CGL-Burgers is spectrally stable by
Definition 3.1.$\qquad\Box$\\

\begin{remark2} We give some certain constants satisfies the
assumption of Theorem 3.1.\\
1). If $m>0$, $r_0=1$, $\vartheta_0=v=0$ and $u=0$ or $u=1$, then we
have
\begin{equation*}
\lambda_1(k)=-k^2m-w_0ki,\;\lambda_2(k)=-k^2-w_0ki,\;\lambda_3(k)=-2-k^2-w_0ki.
\end{equation*}
2). Assume $m>0$, $r_0^2+\vartheta_0^2=1$ and $u=v=0$. Then we have
\begin{equation*}
\lambda_1(k)=-k^2m-w_0ki,\;\lambda_2(k)=-k^2-w_0ki,\;\lambda_3(k)=-(2r_0+k^2+w_0ki).
\end{equation*}
3). If $m>0$, $r_0=u=0$ and $\vartheta_0=1$, then we have
\begin{equation*}
\lambda_1(k)=-k^2m-w_0ki,\;\lambda_{2,3}(k)=-k^2-w_0ki.
\end{equation*}
\end{remark2}

\begin{remark2} Assume that $\kappa=\kappa(r)$, processing  as (3.4), from (3.6).
Then we imply$$ \left|\begin{array}{ccc}

k^2+2r_0^2+[]ki+\lambda&-k^2r_0(c_0+c_1r_0)+2\vartheta_0r_0ki&r_0s_1(r_0)ki\\
c_1k^2+\Gamma&k^2+[]ki+\lambda&\vartheta_0+s_2(r_0)ki \\
2r_0\kappa(r_0)&0&k^2m+w_0ki+\lambda
\end{array}\right|=0,$$
where $\Gamma=c_1\vartheta_0^2+r_0^2v'(r_0)+2r_0v(r_0)$, and
$[]=2\vartheta_0(c_0+c_1r_0)+w_0$. We now consider a simply case as
follows: Letting $r_0=1$, $\vartheta_0=v=0$, $\kappa=m=1$, and
$s_1(r_0)=\frac{1}{8}$. By a long computation,  as $u=c_0+c_1r=0$,
it follows that $\lambda_1(k)=-k^2-w_0ki,$ and
\begin{equation*}
Re\lambda_{2,3}(k)=-(k^2+1)\pm\frac{\sqrt{2}}{4}\sqrt{4+\sqrt{16+k^2}}<0.
\end{equation*}
On the other hand, if $r_0=m=1$, $\vartheta_0=v=0$, $\kappa=1/2$,
and $u=c_0+c_1r=1$, then we have
\begin{equation}
\lambda_{1}(k)=-\frac{1}{3}(2+3k^2+3kw_0i)+\frac{2^{1/3}(4+3ks_1i)}{3\left[a+\sqrt{a^2-4b^3}\right]^{1/3}},
\end{equation}
\begin{equation}\begin{split}
\lambda_{2}(k)=&-\frac{1}{3}(2+3k^2+3kw_0i)-\frac{(1+\sqrt{3}i)(4+3ks_1i)}{3\times2^{2/3}\left[a+\sqrt{a^2-4b^3}\right]^{1/3}}\\&
-\frac{1-\sqrt{3}i}{6\times2^{1/3}}\left(a+\sqrt{a^2-4b^3}\right)^{1/3},
\end{split}\end{equation}
and
 \begin{equation}\begin{split}
\lambda_{3}(k)=&-\frac{1}{3}(2+3k^2+3kw_0i)-\frac{(1-\sqrt{3}i)(4+3ks_1i)}{3\times2^{2/3}\left[a+\sqrt{a^2-4b^3}\right]^{1/3}}\\&
-\frac{1+\sqrt{3}i}{6\times2^{1/3}}\left(a+\sqrt{a^2-4b^3}\right)^{1/3}
\end{split}\end{equation}
with $a=-16+(27k^3s_2-18ks_1)i$ and $b=(4+3ks_1i)$. Therefore, one
can easily check from (3.9), (3.10) and (3.11) that
\begin{equation}
If\quad\left\{\begin{array}{ll}s_1=0, \\
s_2=0,
\end{array}\right.
\quad then
\quad\left\{\begin{array}{ll}\lambda_1(k)=-2-k^2-kw_0i, \\
\lambda_{2,3}(k)=-k^2-kw_0i.
\end{array}\right.
\end{equation}
\begin{equation}
If\quad\left\{\begin{array}{ll}s_1=1, \\
s_2=0,
\end{array}\right.
\quad then
\quad\left\{\begin{array}{ll}\lambda_1(k)=-k^2-kw_0i, \\
\lambda_{2}(k)=-1-k^2-\sqrt{1+ki}-kw_0i,\\
\lambda_{3}(k)=-1-k^2+\sqrt{1+ki}-kw_0i,
\end{array}\right.
\end{equation}
and
\begin{equation}
if\quad\left\{\begin{array}{ll}s_1=-1, \\
s_2=0,
\end{array}\right.
\quad then
\quad\left\{\begin{array}{ll}\lambda_1(k)=-k^2-kw_0i, \\
\lambda_{2}(k)=-1-k^2-\sqrt{1-ki}-kw_0i,\\
\lambda_{3}(k)=-1-k^2+\sqrt{1-ki}-kw_0i.
\end{array}\right.
\end{equation}
By virtue of (3.13) and (3.14), it follows that
$$\lambda_{2,3}(k)=-(1+k^2)\pm\left(\frac{1+\sqrt{1+k^2}}{2}\right)^{1/2}<0.$$
\end{remark2}
\par
Next, we give a usefull lemma which comes from \cite {Kaw, Xin}.

\begin{lemma2} Assume the operator $L$ be simply a lower order perturbation of
the Laplacian, $S(t)$ is an analytic semigroup generated by $L$.
Then the semigroup $S(t)$ satisfies
\begin{equation*}
(i)\quad\|\partial_x^kS(t)u\|_{L^2}\leq
C(1+t)^{-\frac{1}{2}(\frac{1}{2}+k)}\|u\|_{L^1}+Ce^{-\beta
t}\|\partial_x^ku\|_{L^2},
\end{equation*}
\begin{equation*}
(ii)\quad\|\partial_x^kS(t)u\|_{L^2}\leq
C(1+t)^{-\frac{1}{2}(\frac{1}{2}+k)}\|u\|_{L^1}+Ct^{-\frac{1}{2}}e^{-\beta
t}\|\partial_x^{k-1}u\|_{L^2},
\end{equation*}
for some $\beta>0$.
\end{lemma2}
The following lemma will yield the last part of Theorem 3.1.

\begin{lemma2} Assume the initial data $\pi_0$ satisfy $(\|\pi_0\|_{L^1}+\|\pi_0\|_{H^{s+1}})$ is sufficiently small
 for $s>\frac{1}{2}$. Given
$u=0$. If the domain $\mathbb{T}\in \mathbb{R}$ is bounded, then the
solution $\pi$ to Eq.(3.4) satisfies
$$\|\pi\|_{H^{s+1}}\leq C(1+t)^{-\frac{1}{2}(\frac{3}{2}+s)}(\|\pi_0\|_{L^1}+\|\pi_0\|_{H^{s+1}}).$$
\end{lemma2}
\textit{Proof.} Assume $u=c_0+c_1r=0$, by (3.4), we have
\begin{equation}\rho_t=L_{11}\rho-(2\vartheta_0r_0)\phi_x-r_0s_1(r_0)h_x+\psi_1(\rho,h,h_x,\phi_x),
\end{equation}
where $L_{11}=\partial_{x}^2-w_0\partial_x-2r_0^2$, and
$\psi_1(\rho,h,h_x,\phi_x)=-h\rho_x-r_0[\rho^2+\phi_x^2+(s_1(r)-s_1(r_0))h_x]-\rho[2r_0\rho+\rho^2+2\vartheta_0\phi_x+
\phi_x^2+s_1(r)h_x]$.
\par
Since the operator $L_{11}$ is the laplacian with a lower order
perturbation. Thanks to the perturbation theorem in \cite{Pa} to
obtain $L_{11}$ generates an analytic semigroup $S_{11}(t)$, and we
imply from (3.15) that
\begin{equation}
\rho=S_{11}(t)\rho_0+\int_0^tS_{11}(t-\tau)\left[-(2\vartheta_0r_0)\phi_x-r_0s_1(r_0)h_x+\psi_1\right]d\tau.
\end{equation}
Recall that $\|\partial_x^kS_{11}(t)u\|_{L^2}\leq
e^{-r_0^2t}\|\partial_x^k u\|_{L^2}$, one can easily check that
\begin{equation}\left\{\begin{array}{ll}\|S_{11}(t)\rho_0\|_{\dot{H}^k}\leq
e^{-r_0^2t}\|\rho_0\|_{\dot{H}^k}, \\
\int_0^t\|S_{11}(t-\tau)\phi_x\|_{\dot{H}^k}\leq\int_0^te^{-r_0^2(t-\tau)}\|\phi_x\|_{\dot{H}^k}d\tau.
\end{array}\right.
\end{equation}
Applying $\dot{H}^k$ norm to (3.16), in conjunction with (3.17), we
imply
\begin{equation}\|\rho\|_{\dot{H}^k}\leq
e^{-r_0^2t}\|\rho_0\|_{\dot{H}^k}+C\int_0^te^{-r_0^2(t-\tau)}(\|\phi_x\|_{\dot{H}^k}+\|h_x\|_{\dot{H}^k}+
\|\psi_1\|_{\dot{H}^k})d\tau.
\end{equation}
Let $L_{22}=\partial_x^2-w_0\partial_x$. Then the operator $L_{22}$
generates a semigroup $S_{22}(t)$, using (3.4) to yield
\begin{equation}
\phi=S_{22}(t)\phi_0+\int_0^tS_{22}(t-\tau)[-(r_0^2v'(r_0)+2r_0v(r_0))\rho-(s_2(r_0)\partial_x+\vartheta_0)h+\psi_2]d\tau,
\end{equation}
where $\psi_2(\rho,h,\rho_x,h_x,\phi_x)=
-h\phi_x-w_0\vartheta_0-\frac{2\rho_x}{r_0+\rho}(\vartheta_0+\phi_x)-[s_2(r)-s_2(r_0)]h_x-v(r)\rho^2-r_0^2
[v(r)-v'(r_0)\rho]-2r_0\rho[v(r_0)-v(r)].$
\par
In view of lemma 3.2, we obtain
\begin{equation}\begin{split}
\|S_{22}(t)\phi_0\|_{\dot{H}^k}&\leq
C(1+t)^{-\frac{1}{2}(\frac{1}{2}+k)}\|\phi_0\|_{L^1}+Ce^{-\beta
t}\|\phi_0\|_{\dot{H}^k}\\&\leq
C(1+t)^{-\frac{1}{2}(\frac{1}{2}+k)}(\|\phi_0\|_{L^1}+\|\phi_0\|_{\dot{H}^k}),
\end{split}\end{equation}
as $t$ is large enough.\\
Similarly, we can get the following estimation
\begin{equation}\left\{\begin{array}{ll}
\|S_{22}(t-\tau)h_x\|_{\dot{H}^k}\leq
(1+t-\tau)^{-\frac{1}{2}(\frac{1}{2}+k+1)}\|h\|_{L^1}+e^{-\beta(t-\tau)}\|h\|_{\dot{H}^{k+1}}\\
\qquad\leq
C(1+t-\tau)^{-\frac{1}{2}(\frac{1}{2}+k+1)}(\|h\|_{L^1}+\|h\|_{\dot{H}^{k+1}}),\\
\|S_{22}(t-\tau)\psi_2\|_{\dot{H}^{k}}\leq
C(1+t-\tau)^{-\frac{1}{2}(\frac{1}{2}+k)}(\|\psi_2\|_{L^1}+\|\psi_2\|_{\dot{H}^{k}}).
\end{array}\right.
\end{equation}
Therefore, applying $\dot{H}^k$ norm to (3.19), combining (3.20)
with (3.21), it follows that
\begin{equation}\begin{split}
\|\phi\|_{\dot{H}^k}&\leq
(1+t)^{-\frac{1}{2}(\frac{1}{2}+k)}(\|\phi_0\|_{L^1}+\|\phi_0\|_{\dot{H}^k})+\int_0^t
(1+t-\tau)^{-\frac{1}{2}(\frac{1}{2}+k+1)}\\&\times(\|h\|_{L^1}+\|h\|_{\dot{H}^{k+1}})d\tau
+\int_0^t(1+t-\tau)^{-\frac{1}{2}(\frac{1}{2}+k)}\times\\
&(\|\rho\|_{L^1}+\|\rho\|_{\dot{H}^{k}}+\|h\|_{L^1}+\|h\|_{\dot{H}^{k}}+
\|\psi_2\|_{L^1}+\|\psi_2\|_{\dot{H}^{k}})d\tau.
\end{split}\end{equation}
Similarly, we can get the estimation to $h$ as follows
\begin{equation}\begin{split}
\|h\|_{\dot{H}^k}&\leq
(1+t)^{-\frac{1}{2}(\frac{1}{2}+k)}(\|h_0\|_{L^1}+\|h_0\|_{\dot{H}^k})+\int_0^t
(1+t-\tau)^{-\frac{1}{2}(\frac{1}{2}+k)}\\&\times
(\|\psi_3\|_{L^1}+\|\psi_3\|_{\dot{H}^{k}})d\tau
\end{split}\end{equation}
with $\psi_3(\rho,h,\rho_x,h_x)=-hh_x-2\kappa\rho\rho_x$.\\
 Define
$\pi=(\rho,\phi,h)$, if the domain $\mathbb{T}$ is bounded, by
Poincare inequality, then we have $\|\pi\|_{L^2}\leq
C\|\pi\|_{\dot{H}^1}$ and $\|\pi\|_{\dot{H}^k}\sim\|\pi\|_{H^k}$.
Moreover, as $k>\frac{1}{2}$, the space $H^k$ is an algebra. A long
calculation shows that $\psi=(\psi_1,\psi_2,\psi_3)$ satisfies
\begin{equation}
\|\psi\|_{L^1}\leq C\|\psi\|_{\dot{H}^k}\leq
C\|\pi\|_{\dot{H}^{k+1}}^2.
\end{equation}
Adding up (3.18), (3.22) with (3.23), by virtue of (3.24), it
follows that
\begin{equation}\begin{split}\|\pi\|_{\dot{H}^k}&\leq(1+t)^{-\frac{1}{2}(\frac{1}{2}+k)}\|\pi_0\|_{\dot{H}^k}+
\int_0^t(1+t-\tau)^{-\frac{1}{2}(\frac{3}{2}+k)}\|\pi\|_{\dot{H}^{k+1}}d\tau\\&+\int_0^t(1+t-\tau)^{-\frac{1}{2}(\frac{1}{2}+k)}(\|\pi\|
_{\dot{H}^k}+\|\pi\|_{\dot{H}^{k+1}}^2)d\tau.
\end{split}\end{equation}
Analogy to the process of proof to(3.25), one can easily check that
\begin{equation}\begin{split}\|\pi\|_{\dot{H}^{k+1}}&\leq(1+t)^{-\frac{1}{2}(\frac{3}{2}+k)}(\|\pi_0\|_{L^1}+
\|\pi_0\|_{\dot{H}^{k+1}})+\int_0^t(t-\tau)^{-\frac{1}{2}}e^{-\beta(t-\tau)}\|\pi\|_{\dot{H}^{k+1}}^2d\tau\\&+\int_0^t(1+t-\tau)^{-\frac{1}{2}(\frac{3}{2}+k)}
\|\pi\|_{\dot{H}^{k+1}}^2d\tau.
\end{split}\end{equation}
Define $$\mathcal {N}(t)=\sup_{0\leq\tau\leq
t}\{(1+\tau)^{\frac{1}{2}(\frac{3}{2}+k)}\|\pi\|_{\dot{H}^{k+1}}\}.$$
Thanks to (3.26). Using the definition of $\mathcal {N}(t)$, then we
have
\begin{equation}\begin{split}\|\pi\|_{\dot{H}^{k+1}}&\leq(1+t)^{-\frac{1}{2}(\frac{3}{2}+k)}E_0+\mathcal{N}^2(t)
\int_0^t(1+t-\tau)^{-\frac{1}{2}(\frac{3}{2}+k)}
(1+\tau)^{-(\frac{3}{2}+k)}d\tau\\&+\mathcal{N}^2(t)
\int_0^t(t-\tau)^{-\frac{1}{2}}e^{-\beta(t-\tau)}(1+\tau)
^{-(\frac{3}{2}+k)}d\tau.
\end{split}\end{equation}
where $E_0=(\|\pi_0\|_{L^1}+ \|\pi_0\|_{\dot{H}^{k+1}})$.\\
Note that
$$\int_0^t(1+t-\tau)^{-\frac{1}{2}(\frac{3}{2}+k)}
(1+\tau)^{-(\frac{3}{2}+k)}d\tau\leq
C(1+t)^{-\frac{1}{2}(\frac{3}{2}+k)},$$
$$\int_0^t(t-\tau)^{-\frac{1}{2}}e^{-\beta(t-\tau)}(1+\tau)
^{-(\frac{3}{2}+k)}d\tau\leq C(1+t)^{-(\frac{3}{2}+k)}.$$ Therefore
\begin{equation}
\mathcal{N}(t)\leq CE_0+C_1\mathcal{N}^2(t).
\end{equation}
In view of the assumption stated in the lemma, we imply that
$\mathcal{N}(t)\leq CE_0$. This completes the proof the lemma.$\qquad\Box$\\
\begin{remark2} Assume that the initial $\pi_0=(\rho_0,\phi_0,h_0)\in
(H^{s+1})^3,\;s>\frac{1}{2}$. Then Eq.(3.4) becomes $\pi_t=\mathcal
{L}\pi +\psi$, and the operator $\mathcal {L}$ generates a semigroup
$S(t)$. Due to $\|\partial_x^sS(t)\pi\|_{L^2}\leq
Ct^{-\frac{s}{2}}\|\pi\|_{L^2}$,
 $H^s$ is an algebra, and the higher order term $\psi$ is continuous in $H^s$, by the semigroup theory \cite{Pa},
 there exists a positive $T>0$, such that the solution $\pi\in \mathcal {C}([0,T[;H^{s+1})\cap
 \mathcal {C}^1([0,[;H^s)$. This ensures the existence of solution to Eq.(3.4) in Lemma 3.3.
\end{remark2}

\section{Nonlinear Instability}
 In this section, we present the instability of the equilibrium solution of the CGL--Burgers equations by the principle of linearized
instability in \cite{Lun}. We consider the following equation
\begin{equation}\left\{\begin{array}{ll} \pi_t=\mathcal
{L}\pi+\psi,&(t,x)\in\mathbb{R}^{+}\times\mathbb{T},\\
\pi|_{t=0}=\pi_0=(\rho_0,\phi_0,h_0),&x\in\mathbb{T},
\end{array}\right.\end{equation}
where the operator $\mathcal {L}=\mathcal {A}\partial_x^2+\mathcal
{B}\partial_x+\mathcal {C}$, the functions $\pi=(\rho,\phi,h)$ and
$\psi=(\psi_1,\psi_2,\psi_3)$,
\begin{equation}\left\{\begin{array}{ll}\psi_1(\rho,h,h_x,\phi_x)=-2\vartheta_0c_1\rho\rho_x-2(c_0+c_1r)\phi\rho_x-
(c_0+c_1r_0)\rho\phi_{xx}-h\rho_x\\
\quad-r_0[\rho^2+\phi_x^2+(s_1(r)-s_1(r_0))h_x]
-\rho[2r_0\rho+\rho^2+2\vartheta_0\phi_x+\phi_x^2+s_1(r)h_x],\\
\psi_2(\rho,h,\rho_x,h_x,\phi_x)=
-h\phi_x-w_0\vartheta_0-(c_0+c_1r_0)(\vartheta_0^2+\phi_x^2)+\frac{c_0}{r_0+\rho}\rho_{xx}\\
\quad-c_1(2\vartheta_0\phi_x+\phi_x^2)-\frac{2\rho_x}{r_0+\rho}(\vartheta_0+\phi_x)-v(r)\rho^2-[s_2(r)-s_2(r_0)]h_x\\
\quad-r_0^2[v(r)-v'(r_0)\rho]-2r_0\rho[v(r_0)-v(r)],\\
\psi_3(\rho,h,\rho_x,h_x)=-hh_x-2\kappa\rho\rho_x.
\end{array}\right.\end{equation}
We shall prove an instability result under the assumption
\begin{equation}\left\{\begin{array}{ll}\sigma_+(\mathcal {L})=\sigma(\mathcal
{L})\cap\{\lambda\in\mathbb{C}:Re\lambda>0\}\neq\emptyset,\\
\inf\{Re\lambda:\lambda\in\sigma_+(\mathcal {L})\}=\omega_+>0.
\end{array}\right.\end{equation}
\begin{theorem3}
 Under the condition (4.3). There exists a bounded domain
 $\mathbb{T}\in\mathbb{R}$ such that the equilibrium solution of the
coupled CGL--Burgers equations is unstable.
\end{theorem3}
By virtue of Theorem 9.1.3 of page 344 in \cite{Lun}, we will break
the result into two lemmas.
\begin{lemma3} The operator $\mathcal {L}=\mathcal {A}\partial_x^2+\mathcal
{B}\partial_x+\mathcal {C}$ is sectorial and the graph morn of
$\mathcal {L}$ is equivalent to the morn of $D(\mathcal {L})$.
Moreover, there exists domain $\mathbb{T}$ and a positive number
$\delta_0$ such that
\begin{equation}\left\{\begin{array}{ll}\sigma_+(\mathcal
{L})=\sigma(\mathcal
{L})\cap\{\lambda\in\mathbb{C}:Re\lambda>0\}\neq\emptyset,\\
\inf\{Re\lambda:\lambda\in\sigma_+(\mathcal {L})\}=\delta_0>0.
\end{array}\right.\end{equation}
\end{lemma3}
\textit{Proof.} From (3.7) and (3.8), we have $Re\lambda_{2}(k)<0$.
Case 1: If $k^2m>0$, $2(r_0^2+k^2)^2< a$, or
$4(r_0^2+k^2)^4-4a(r_0^2+k^2)^2<b^2$, then we deduce that
$$Re\lambda_1(k)<0,\quad\text{and} Re\lambda_{3}(k)>0.$$
Therefore, there exists $\Omega\in\mathbb{R}$ such that
$$\inf\{Re\lambda:\lambda\in\sigma_+(\mathcal {L})\}>0.$$
 Case 2: If $m<0$,  then it follows that
$Re\lambda_1(k)>0$. Therefore, there also exists
$\Omega\in\mathbb{R}$ such that
$$\inf\{Re\lambda:\lambda\in\sigma_+(\mathcal {L})\}>0.$$
Since the Laplacian is a sectorial operator, and the operator
$\mathcal {L}$ is a lower order perturbation of the Laplacian, by
the perturbation theorem, it is easy to show that $\mathcal {L}$ is
a sectorial operator and the graph of $\mathcal {L}$
is equivalent to the norm of $D(\mathcal {L})$.$\qquad\Box$\\
\par
In order to obtain Theorem 4.1, we only to prove the following lemma
by virtue of Theorem 9.1.3 of page 344 in \cite{Lun}.
\begin{lemma3} Suppose $\mathbb{D}$ be a neighborhood of the origin
in $D(\mathcal {L})$. Then there exists constants $r_0,\vartheta_0$
and the functions $u,v$ such that the map
$\psi:\mathbb{D}\longrightarrow L^2$ is a $\mathcal {C}^1$ function
with locally Lipschitz continuous derivative and satisfies
$$\psi(0)=0,
\qquad \psi'(0)=0,$$ where $\psi'(0)$ is the Fr\'{e}chet derivative
of $\psi(\pi)$ at origin.
\end{lemma3}
\textit{Proof.} First, by virtue of definition of $\psi(\pi)$, we
have $\psi(0)=(0,-w_0\vartheta_0-(c_0+c_1r_0)\vartheta_0^2,0)$.
Therefore, we deduce $w_0\vartheta_0+(c_0+c_1r_0)\vartheta_0^2=0$.
If $D_1,D_2,D_3$ denote the Fr\'{e}chet derivative of $\rho,\phi,h$,
respectively. Then we obtain
\begin{equation}\left\{\begin{array}{ll}D_1\psi_1=-2\vartheta_0c_1\rho_x-2\phi\rho_x
-(c_0+c_1r_0)\phi_{xx}-r_0[2\rho+s_1'(r)h_x]\\\quad
-[4r_0\rho+3\rho^2+2\vartheta_0\phi_x+\phi_x^2+s_1(r)h_x+\rho
s_1'(r)h_x]-[2\vartheta_0c_1\rho\\\quad+2(c_0+c_1r)\phi+h]\partial_x\\
\;\;\qquad=H_{11}+J_{11}\partial_x,\\
D_2\psi_1=-2(c_0+c_1r)\rho_x-(2r_0\phi_x+2\vartheta_0\rho+2\rho\phi_x)\partial_x-(c_0+c_1r_0)\rho\partial_x^2\\
\;\;\qquad=H_{12}+J_{12}\partial_x+K_{12}\partial_x^2,\\
 D_3\psi_1=-\rho_x+[r_0(s_1(r_0)-s_1(r))-\rho
s_1(r)]\partial_x\\
\;\;\qquad=-\rho_x+J_{13}\partial_x.
\end{array}\right.\end{equation}
\begin{equation}\left\{\begin{array}{ll}
D_1
\psi_2=-\frac{c_0\rho_{xx}}{(r_0+\rho)^2}-c_1(2\vartheta_0\phi_x+\phi_x^2)+2\frac{(\vartheta_0+\phi_x)\rho_x}
{(r_0+\rho)^2}-v'(r)\rho^2\\\qquad-2v(r)\rho-r_0^2[v'(r)-v'(r_0)]-2r_0[v(r_0)-v(r)]+2r_0v'(r)\rho
\\\qquad-2\frac{(\vartheta_0+\phi_x)\partial_x}{r_0+\rho}+\frac{c_0\partial_x^2}{r_0+\rho}\\
\;\;\qquad=H_{21}+J_{21}\partial_x+K_{21}\partial_x^2,\\
D_2\psi_2=-\left[h+2(c_0+c_1r_0)\phi_x+2c_1\rho(\vartheta_0+\phi_x)+\frac{2\rho_x}{r_0+\rho}\right]\partial_x\\
\qquad\;\;=J_{22}\partial_x,\\
D_3\psi_2=-\phi_x+[s_2(r_0)-s_2(r)]\partial_x\\
\;\;\qquad=-\phi_x+J_{23}\partial_x.\\
\end{array}\right.\end{equation}

\begin{equation}\left\{\begin{array}{ll}
D_1\psi_3=-2\kappa\rho_x-2\kappa\rho\partial_x,\\
D_2\psi_3=0,\\
D_3\psi_3=-h_x-h\partial_x.
\end{array}\right.\end{equation}
Therefore,
$$\psi'(\pi)|_{\pi=0}=\left(\begin{array}{ccc}

D_1\psi_1&D_2\psi_1&D_3\psi_1 \\
D_1\psi_2&D_2\psi_2&D_3\psi_2 \\
D_1\psi_3&D_2\psi_3&D_3\psi_3
\end{array}\right)_{|_{\pi=0}}=
\left(\begin{array}{ccc}

0&0&0 \\
-2\frac{\vartheta_0}{r_0}\partial_x+\frac{c_0}{r_0}\partial_x^2&0&0 \\
0&0&0
\end{array}\right).$$
Then, we have $\vartheta_0=c_0=0$ and $v=0,r_0^2=1$.
\par
In view of (4.5), (4.6) and (4.7), we can write $\psi'(\pi)$ as
follows
$$\psi'(\pi)=\left(\begin{array}{ccc}
H_{11}&H_{12}&-\rho_x \\
H_{21}&0&-\phi_x\\
-2\kappa\rho_x&0&-h_x
\end{array}\right)+\left(\begin{array}{ccc}

J_{11}&J_{12}&J_{13} \\
J_{21}&J_{22}&J_{23}\\
-2\kappa\rho&0&-h
\end{array}\right)\frac{\partial}{\partial x}$$$$+
\left(\begin{array}{ccc}
0&K_{12}&0 \\
K_{21}&0&0\\
0&0&0
\end{array}\right)\frac{\partial^2}{\partial x^2}=H(\pi)+J(\pi)\frac{\partial}{\partial x}+K(\pi)
\frac{\partial^2}{\partial x^2}$$
Let $\pi_1=(\rho_1,\phi_1,h_1)\in\mathbb{D}$,
$\pi_2=(\rho_2,\phi_2,h_2)\in\mathbb{D}$, for any
$R=(f,g,\eta)\in\mathbb{D}$, we imply that
\begin{equation}\begin{split}
\|[\psi'(\pi_1)-&\psi'(\pi_2)]R\|\leq\|[H(\pi_1)-H(\pi_2)]R\|+\\&\|[J(\pi_1)-J(\pi_2)]R_x\|+
\|[K(\pi_1)-K(\pi_2)]R_{xx}\|.
\end{split}\end{equation}
Because of
\begin{equation}\begin{split}\|[H(\pi_1)&-H(\pi_2)]R\|=\sum_{i=1}^3(\|[H_{i1}(\pi_1)-H_{i1}(\pi_2)]f\|+
\\&\|[H_{i2}(\pi_1)-H_{i2}(\pi_2)]g\|+\|[H_{i3}(\pi_1)-H_{i3}(\pi_2)]\eta\|)\\
&\qquad\leq
C(\|\pi\|_{H^{1+s}})\|\pi_1-\pi_2\|_{H^1}(\|f\|_{L^2}+\|g\|_{L^2}+\|\eta|_{L^2}),
\end{split}\end{equation}
\begin{equation}\begin{split}\|[J(\pi_1)&-J(\pi_2)]R_x\|=\sum_{i=1}^3(\|[J_{i1}(\pi_1)-J_{i1}(\pi_2)]f_x\|+
\\&\|[J_{i2}(\pi_1)-J_{i2}(\pi_2)]g_x\|+\|[J_{i3}(\pi_1)-J_{i3}(\pi_2)]\eta_x\|)\\
&\qquad\leq
C(\|\pi\|_{H^{1+s}})\|\pi_1-\pi_2\|_{H^1}(\|f_x\|_{L^2}+\|g_x\|_{L^2}+\|\eta_x|_{L^2}),
\end{split}\end{equation}
and
\begin{equation}\begin{split}&\|[K(\pi_1)-K(\pi_2)]R_{xx}\|=\sum_{i=1}^3(\|[K_{i1}(\pi_1)-K_{i1}(\pi_2)]f_{xx}\|+
\\&\qquad\qquad\|[K_{i2}(\pi_1)-K_{i2}(\pi_2)]g_{xx}\|+\|[K_{i3}(\pi_1)-K_{i3}(\pi_2)]\eta_{xx}\|)\\
&\leq C(\|\pi\|_{H^{s}})\|\rho_1-\rho_2\|_{H^1}\|g_{x}\|_{L^2},
\end{split}\end{equation}
where we used the integrating by parts and Sobolev's embedding
theorem.\\

 Combining (4.8), (4.9), (4.10) with (3.11), it follows that
$$\|\psi'(\pi_1)-\psi'(\pi_2)\|_{H^{-1}}\leq
C(\|\pi\|_{H^{1+s}})\|\pi_1-\pi_2\|_{H^1}.$$
This completes the proof Lemma 4.2. $\qquad\Box$\\

 \textbf{Acknowledgments}

This work was partially supported by CPSF (Grant No.: 2012M520007).
The authors thank the references for their valuable comments and
constructive suggestions.\\

\end{document}